\documentclass[a4paper,two-side,11pt]{article}
\usepackage[latin1]{inputenc}
\usepackage[T1]{fontenc}
\usepackage{amsmath}
\usepackage{amsthm}
\usepackage{latexsym}
\usepackage{amssymb}
\usepackage[all]{xy}
\usepackage{calc}
\usepackage{multicol}
\newtheorem{thm}{Theorem}[section]
\newtheorem{prop}{Proposition}[section]

\newtheorem{defi}{Definition}[section]

\newtheorem{lem}{Lemma}[section]
\newtheorem{cor}{Corollary}[section]

\newtheorem{etape}{Step}
\newtheorem{rmq}{Remark}[section]

\newcommand{\R}{\mathbb{R}}
\numberwithin{equation}{section}
\newcommand{\N}{\mathbb{N}}

\newcounter{exercice}

\author{Jean-baptiste Castéras}
\title{A mean field type flow.}
\begin{document}
\date{}
\maketitle
 \begin{center}
 \small{LMBA\\
 Universit\'e de Bretagne Occidentale\\
 6 av. Victor Le Gorgeu -- CS 93837\\
 29238 Brest Cedex\\
 France\\
 E-mail:  jean-baptiste.casteras@univ-brest.fr}
 \end{center}
\medskip
 \begin{abstract}
We consider a gradient flow related to the mean field type equation. First, we show that this flow exists for all time. Next, we prove a compactness result for this flow allowing us to get, under suitable hypothesis on its energy, the convergence of the flow to a solution of the mean field type equation. We also get a divergence result if the energy of the initial data is largely negative.
 \end{abstract}
 \vspace{10pt}
 \begin{center} \textit{Key words} : Mean field equation, Blow-up analysis, Geometric flow.
 \end{center}
 \begin{center} \textbf{AMS subjet classification} : 35B33, 35J20, 53C44, 58E20.
 \end{center}
\subsection{Introduction}
Let $(M,g)$ be a compact riemannian surface without boundary, we will study an evolution problem associated to a mean field type equation :
\begin{equation}
\label{eqell}
-\Delta v +Q=\rho \frac{e^v}{\int_M e^v dV},
\end{equation}

\noindent where $\rho$ is a real number, $Q\in C^\infty (M)$ is a given function such that $\displaystyle\int_M Q dV=\rho$ and $\Delta$ is the Laplacian with respect to the metric $g$. Equation (\ref{eqell}) is equivalent to the mean field equation :
\begin{equation}
\label{E:CM}
  -\bigtriangleup u + \rho  \left( \dfrac{-f  e^{u}}{\int_M f e^{u}dV} +\dfrac{1}{|M|}\right)=0 ,
 \end{equation}
where $|M|$ stands for the volume of $M$ with respect to the metric $g$ and $f\in C^\infty(M)$ is a function supposed positive.
Indeed, if $v$ is a solution of \eqref{eqell}, by setting $v=u+\log f$, we recover that $u$ is a solution of \eqref{E:CM} with $Q=\dfrac{\rho}{|M|}+\Delta \log f$. The mean field equation appears in conformal geometry but also in statistic mechanic from Onsager's vortex model for turbulent Euler flows. More precisely, in this setting, the solution $u$ of the mean field equation is  the stream function in the infinite vortex limit (see \cite{MR1145596}). It also arises in the abelian Chern-Simons-Higgs model (see for example \cite{MR1324400}, \cite{MR1955272}, \cite{MR1400816}, \cite{MR1838682}).\newline

Equation (\ref{E:CM}) has a variational structure and its solutions can be found as the critical points of the following functional :
 \begin{equation}
 \label{enerE:CM}
 I_\rho (u)=\frac{1}{2}\int_M |\nabla u|^2 dV +\frac{\rho}{|M|}\int_M udV -\rho \log \left( \int_M fe^u dV\right),\ u\in H^1(M).
 \end{equation}
When $\rho<8\pi$, from the Moser-Trudinger's inequality, one can easily prove that the functional $I_\rho$ is bounded from below and coercive ; thus one can find solution of (\ref{E:CM}) by minimizing $I_\rho$. Existence of solutions becomes more delicate if $\rho\geq 8\pi$. When $\rho=8\pi$, $I_\rho$ admits a lower bound and isn't anymore coercive and when $\rho>8\pi$, $I_\rho$ is  unbounded from below and above. Existence of solutions of equation \eqref{eqell} has been widely studied this last decades. Many partial existence results have been obtained if $\rho\neq 8k\pi,\ k\in \N^\ast$, and according to the Euler characteristic of $M$ (see for example \cite{MR1132783}, \cite{MR2001443}, \cite{MR1712560}, \cite{MR1673972}, \cite{MR1322618}, \cite{MR2483132}, \cite{MR1619043}). Finally, when $\rho\neq 8k\pi$, $k\in \N^\ast$, Djadli \cite{MR2409366} has generalized the previous results establishing the existence of solution for all surfaces $M$ by studying the topology of sublevels $\left\{I_\rho \leq -C\right\}$ to achieve a min-max scheme (already introduced in Djadli-Malchioldi \cite{MR2456884}).\newline

In this paper, we consider the evolution problem associated to equation (\ref{eqell}), that is the following equation
\begin{equation}
\label{E:flot}
\left\{\begin{array}{ll}
 \frac{ \partial}{\partial t}e^{v} =\bigtriangleup v - Q+ \rho \dfrac{  e^{v}}{\int_M  e^{v}dV} \\
 v(x,0)=v_{0}(x),
\end{array}
\right.
\end{equation}
where $v_0 \in C^{2+\alpha}(M)$, $\alpha\in (0,1)$, is the initial data and $Q\in C^\infty (M)$ is a given function such that $\displaystyle\int_M Q dV=\rho$. It is a gradient flow with respect to the following functional which will be called energy :
\begin{equation}
\label{E:fonc}
J_\rho(v)= \frac{1}{2}\int_M|\nabla v|^2dV+\int_M  Qv(t) dV- \rho \ln\left(\int_M  e^{v}dV\right),\ v\in H^1(M). 
\end{equation}
This functional is unbounded from below (except in the case $\rho < 8\pi$) and above. The interest of this flow is that it satisfies some important geometrical properties useful for its convergence (see in particular estimate \eqref{mp4} of Section $3$). When $Q$ is a constant equal to the scalar curvature of $M$ with respect to the metric $g$, the flow \eqref{E:flot} (normalized) has been studied by Struwe \cite{MR1991140} (we notice that in this case $\rho=\displaystyle\int_M QdV\leq 8\pi$). For others curvature flows, we refer to \cite{MR2064965}, \cite{MR1999924}, \cite{MR2168505}, \cite{MR2254307}, \cite{MR954419}, \cite{MR2559723}, \cite{MR2217518}, \cite{MR2011332}, \cite{MR2137948} and the references therein.\newline
We begin by studying the global existence of the flow (\ref{E:flot}). We prove the following theorem :
\begin{thm}
\label{exglob} 
For all initial data $ v_0 \in C^{2+\alpha}(M)$, $\alpha\in (0,1)$, all $\rho\in \R$ and all function $Q\in C^\infty (M)$ such that $\displaystyle\int_M Q dV=\rho$, there exists a unique global solution $v \in C_{loc}^{2+\alpha,1+\frac{\alpha}{2}}\left(M\times [0, + \infty)\right)$ of \eqref{E:flot}.
\end{thm}
Next, we investigate the convergence of the flow. We denote $v(t): M\rightarrow \R$ the function defined by $v(t)(x)=v(x,t)$. We will show that if the energy $J_\rho (v(t))$ of the global solution is bounded from below uniformly in time (when $\rho >8\pi$) then $v(t)$ converges when $t\rightarrow +\infty$ to a function $v_\infty$ solution of equation \eqref{eqell}. More precisely, we have :
\begin{thm}
\label{conv1}
Let $v(t)$ be the solution of \eqref{E:flot}.\newline
(i) If $\rho <8\pi$, then $v(t)$ converges in $H^2(M)$ to a function $v_\infty \in C^\infty (M)$ solution of
\begin{equation*}
-\Delta v_\infty +Q=\rho \frac{e^{v_\infty}}{\int_M e^{v_\infty} dV}.
\end{equation*}
(ii) If $\rho> 8\pi$, $\rho \neq 8k\pi$, $k\in \N^\ast$ and if there exists a constant $C>0$ not depending on $t$ such that, for all $t\geq 0$,
\begin{equation}
\label{29juin2012e2}
J_\rho (v(t))\geq -C,
\end{equation}
then $v(t)$ converges in $H^2(M)$ to a function $v_\infty \in C^\infty (M)$ solution of
\begin{equation*}
-\Delta v_\infty +Q=\rho \frac{e^{v_\infty}}{\int_M e^{v_\infty} dV}.
\end{equation*}
\end{thm} 
Moreover, we will prove that there exist initial data $v_0\in C^{\infty}(M)$ such that the energy of the global solution $v(t)$ of the flow, with $v(0)(x)=v_0(x)$ for all $x\in M$, stays uniformly bounded from below hence, thanks to Theorem \ref{conv1}, such that the flow converges. 
\begin{thm}
\label{conv2}
Let $\rho\neq8k\pi$, $k\in \N^\ast$. There exist initial data $v_0\in C^{\infty}(M)$ such that the global solution $v(t)$ of (\ref{E:flot}) with $v(x,0)=v_0(x)$, for all $x\in M$, satisfies \eqref{29juin2012e2} i.e. such that the global solution $v(t)$ of (\ref{E:flot}) with $v(x,0)=v_0(x)$, for all $x\in M$, converges in $H^2(M)$ to a function $v_\infty \in C^\infty (M)$ solution of the equation \eqref{eqell}
 \begin{equation*}
-\Delta v_\infty +Q=\rho \frac{e^{v_\infty}}{\int_M e^{v_\infty} dV}.
\end{equation*}
\end{thm}

Finally, we will show that if the energy of the initial data $v_0$ of (\ref{E:flot}) is largely negative then the flow diverges when $t\rightarrow +\infty$.
\begin{thm}
\label{conv3}
Let $\rho \in ( 8k\pi, 8(k+1)\pi)$ ($k\geq 1$). Then, there exists a constant $C>0$ depending on $M,Q$ and $\rho$ such that for all $ v_{0}\in C^{2+\alpha}(M)$ satisfying $J_\rho(v_{0})\leq -C$, then the global solution $v(t)$ of (\ref{E:flot}) with $v(x,0)=v_0(x)$, for all $x\in M$, satisfies $J_\rho (v(t))\underset{n\rightarrow +\infty}{\longrightarrow } -\infty$.
\end{thm}
In order to prove the previous convergence results, we study the compactness property of solutions $(v_n)_n \subseteq H^2 (M)$ of the following perturbed elliptic mean field type equation :
\begin{equation}
\label{regjanveqini}
-\Delta v_n =Q_0 +h_n e^{v_n}+\rho e^{v_n},
\end{equation}
where $\rho>0$, $Q_0 \in C^0 (M)$ is a given function and $(h_n)_n \subseteq C^0(M)$. In fact, in what follows, the sequence $(v_n)_n\subseteq H^2(M)$ is said compact if it is uniformly bounded in $H^2(M)$. The term $h_n$ corresponds to the parabolic term of \eqref{E:flot}. We also assume that there exists a constant $C>0$ not depending on $n$ such that
\begin{equation}
\label{regdeceq1}
\left\{\begin{array}{ll}
(i)\displaystyle\lim_{n\rightarrow +\infty}\int_M h_n^2 e^{v_n}dV = 0, \\
(ii)\ h_n(x) e^{v_n(x)}+\rho e^{v_n(x)}\geq -C,  \forall x\in M,\ \forall n\geq 0.
\end{array}
\right.
\end{equation}
We will see that these conditions are satisfied by the solution of the flow \eqref{E:flot}. We obtain the following compactness result : 
\begin{thm}
\label{regdecthm3}
Let $(v_n)_n \subseteq H^2(M)$ be a sequence of solutions of \eqref{regjanveqini} such that $\displaystyle\int_M e^{v_n}dV=1$, $\forall n\geq 0$, and satisfying  conditions \eqref{regdeceq1}. Then, we have the following alternative :\\

- Either, there exists a constant $C$ not depending on $n$ such that
$$\left\|v_n \right\|_{H^2 (M)}\leq C, \forall n\in \N,$$
 
- or, up to a subsequence, there exist $k$ sequences of points $(x_n^1)_{n},\ldots ,(x_n^k)_{n}$ with $1\leq k\leq [\dfrac{\rho}{8\pi}]$ (where $[\dfrac{\rho}{8\pi}]$ stands for the integer part of $\dfrac{\rho}{8\pi}$), and $k$ sequences of real positive numbers $(R_n^1)_n,\ldots ,(R_n^k)_n$ satisfying the following properties\newline
(i)
\begin{equation}
\label{regdecprop3.1.bister}
 \lim_{n\rightarrow +\infty} \int_{B_{2R_n^i}(x_n^i)} e^{v_n}dV=\dfrac{8\pi}{\rho},\ \forall i\in \{1,\ldots , k\}.
 \end{equation}
where $B_{2R_n^i}(x_n^i)$ stands for the geodesic ball of center $x_n^i$ and radius $2R_n^i$.\\ 
(ii)
\begin{eqnarray}
\label{regdecprop3.2.bister} 
&& \dfrac{R_n^i}{|x_n^i-x_n^j|}\underset{n\rightarrow +\infty}{\longrightarrow}  0,\ \forall i\neq j\in \{1,\ldots ,k\}\ if\ k\neq 1,\nonumber\\
 && R_n^1\underset{n\rightarrow +\infty}{\longrightarrow}  0,\ if\ k=1,
 \end{eqnarray}
 where $|x_n^i-x_n^j|$ stands for the geodesic distance from $x_n^i$ to $x_n^j$ with respect to the metric $g$.\\
(iii)
\begin{equation*}
\label{26maieqsup1}
\displaystyle\lim_{n\rightarrow +\infty} \int_{M\backslash \bigcup_{i=1}^k B_{2R_n^i}(x_n^i)}e^{v_n}dV=0.
\end{equation*}
\end{thm}

The compactness of solutions of \eqref{regjanveqini} when $\rho \neq 8k\pi$, $k\in \N^\ast$, is an immediate consequence of Theorem \ref{regdecthm3}.
\begin{cor}
\label{regdecthm2}
Let $(v_n)_n\subseteq H^{2}(M)$ be a sequence of solutions of (\ref{regjanveqini}) such that $\displaystyle\int_M e^{v_n}dV=1$, $\forall n\geq 0$, and satisfying \eqref{regdeceq1}. If $\rho\neq 8k\pi$, $k\in \N^\ast$, then there exists a constant $C$ not depending on $n$ such that
$$\left\|v_n\right\|_{H^2(M)}\leq C.$$ 
\end{cor} 
Corollary \ref{regdecthm2} generalizes the result of Li \cite{MR1673972} where the function $h_n$ is equal to zero. Due to the term $h_n e^{v_n}$, the methods used in Li \cite{MR1673972} are not adapted anymore, therefore we will rely on integral estimates in the spirit of Malchiodi \cite{MR2248155}.\\ 

The paper is organised as follows : in Section $1$, we prove the global existence of solution of \eqref{E:flot}. We also show the continuity of the flow with respect to its initial data. In Section $2$, we study the compactness property of $(v_n)_n \subseteq H^2(M)$ solution of \eqref{regjanveqini}. We first show that we have the following alternative : either $(v_n)_n$ is compact or some blow-up phenomena for $(v_n)_n$ appear. Next we establish the asymptotic profil of $(v_n)_n$ near a blow-up. Then we prove an integral Harnack type inequality which will allow us to obtain estimates $(i),(ii)$ and $(iii)$ of Theorem \ref{regdecthm3}. Corollary \ref{regdecthm2} is an immediate consequence of Theorem \ref{regdecthm3} since if we assume that we are in the second alternative of Theorem \ref{regdecthm3} then it follows that
$$1=\int_M e^{v_n}dV \underset{n\rightarrow +\infty}{\longrightarrow} \dfrac{8k\pi}{\rho}.$$
Therefore, we deduce that if $\rho\neq 8k\pi$, $k\in \N^\ast$, the concentration phenomena can't appear and, thereby, that the sequence $(v_n)_n\subseteq H^2(M)$ is compact.\newline 
In Section $3$, we study the convergence of the flow \eqref{E:flot}. We begin by proving Theorem \ref{conv1}. We first show that the global solution $v(t)$ of \eqref{E:flot} is uniformly (with respect to $t$) bounded in $H^1(M)$ when $v(t)$ satisfies condition \eqref{29juin2012e2}. The proof involves the compactness result obtained in Corollary \ref{regdecthm2}. We point out that, when $\rho <8\pi$, condition \eqref{29juin2012e2} is always satisfied. Then, we show that the parabolic term of \eqref{E:flot}, $\dfrac{\partial v(t)}{\partial t}$, tends to $0$ when $t\rightarrow +\infty$ in $L^2(M)$ norm with respect to the metric $g_1(t)=e^{v(t)}g$. This implies that $v(t)$ is uniformly bounded in $H^2(M)$. Next we prove Theorem \ref{conv2}, i.e. there exists initial data in $C^\infty (M)$ such that condition \eqref{29juin2012e2} is satisfied. Our proof is based on the study of the topology of the level set
$$\left\{v\in X\ :\ J_\rho (v) \leq -L \right\},$$
where $X$ is the space of $C^\infty (M)$ functions endowed with the $C^{2+\alpha}(M)$ norm, $\alpha\in (0,1)$. The end of Section $3$ is devoted to the proof of Theorem \ref{conv3}.
\vspace{10pt}
\section{Global existence of the flow.}
We begin by noticing that, since the flow is parabolic, standard methods (see for example \cite{MR0181836}) provide short time existence. 
Thus, there exists $T_1>0$ such that $v\in C^{2+\alpha,1+\frac{\alpha}{2}}(M\times [0,T_1])$ is a solution of \eqref{E:flot}. We will give two basic properties of the flow : the conservation of the volume of $M$ endowed with the metric $g_1(t)=e^{v(t)}g$ and the decreasing, along the flow, of the functional $J_\rho(v(t))$.

\begin{prop}
(i) For all $ t\in [0,T_1]$, we have
\begin{equation}
\label{volconst}
\int_M e^{v(t)}dV=\int_M e^{v_0}dV.
\end{equation}
(ii) If  $0\leq t_0\leq t_1\leq T_1$, we have :
\begin{equation}
\label{eg2.5}
 J_\rho(v(t_1))\leqslant J_\rho(v(t_0)).
 \end{equation}
\end{prop}
\begin{proof}
To see that \eqref{volconst} holds, it is sufficient to integrate (\ref{E:flot}) on $M$. Differentiating $J_\rho (v(t))$ with respect to $t$ and integrating by parts, one finds, for all $t\in [0,T_1]$,
\begin{equation}
\label{eg2.5'}
\dfrac{\partial}{\partial t}J_\rho(v(t))= -\int_M \left( \dfrac{\partial v(t)}{\partial t} \right)^2 e^{v(t)}dV \leqslant 0.  \hspace{2cm} 
\end{equation}
This implies \eqref{eg2.5}.
\end{proof}

\subsection {Proof of Theorem \ref{exglob}.}
To prove the global existence of the flow, we set 
$$T=\sup\left\{\overline{T}>0 :\exists v \in C^{2+\alpha, 1+\frac{\alpha}{2}}(M\times [0,\overline{T}])\ solution\ of\ (\ref{E:flot})\right\},$$ 
and suppose that $T<+\infty$. From the definition of $T$, we have $v\in C_{loc}^{2+\alpha, 1+\frac{\alpha}{2}}(M\times [0,T))$. We will show that there exists a constant $C_T>0$ depending on $T$, $M$, $Q$, $\rho$ and $\left\|v_0\right\|_{C^{2+\alpha}(M)}$ such that
\begin{equation}
\label{7dec2012e1}
\left\|v\right\|_{C^{2+\alpha, 1+\frac{\alpha}{2}}(M\times [0,T))}\leq C_T.
\end{equation}
This estimate allows us to extend $v$ beyond $T$, contradicting the definition of $T$.\\

In the following, $C$ will denote constants depending on $M$, $Q$, $\rho$ and $\left\|v_0\right\|_{C^{2+\alpha}(M)}$ and $C_T$ the ones depending on $M$, $Q$, $\rho$, $\left\|v_0\right\|_{C^{2+\alpha}(M)}$ and $T$.
\begin{prop}
\label{etape4}
For all $\rho \in \R$, there exists a constant $C_T$ such that
\begin{equation}
\label{etape4e}
\left\|v(t)\right\|_{H^1(M)}\leq C_T,\ \forall t\in [0,T).
\end{equation}
Moreover, if $\rho < 8\pi$, there exists a constant $C$ not depending on $T$ such that
\begin{equation}
\label{etape4e1}
\left\|v(t)\right\|_{H^1(M)}\leq C,\ \forall t\in [0,T).
\end{equation}
\end{prop}
\begin{proof}
We will decompose the proof into three steps.
\newcounter {EqNo}
 \setcounter{EqNo}{0}
\newcommand{\Numeq}{\refstepcounter{EqNo}  \hfill ( \thesection.\theEqNo)\\ }
\newcommand{\refeq}[1]{(\thesection.\ref{#1})}
\noindent 
\begin{etape}
\label{etape1}
\label{vborn}
Let $\rho\geq 8\pi$. There exists a constant $C_T$ such that
\begin{equation}
\label{vborne1}
 v(x,t) \leqslant C_T,\ \forall x \in M ,\ \forall t\in [0,T).
\end{equation}
\end{etape}

\noindent\textit{Proof of Step 1.} We define $v_{\max}(t)=\underset{x\in M}{\max}\ v\left(x,t\right)$. By the maximum principle and since $v$ satisfies equation \eqref{E:flot}, we find
 \begin{eqnarray*}
 \label{fardoun15}
 \dfrac{\partial }{\partial t}e^{v_{\max}(t)}  & \leq &\frac{\rho}{\int_{M}e^{v_0}dV}\left( \left\|Q\right\|_{L^\infty (M)} \frac{\int_{M}e^{v_0}dV}{\rho}+e^{v_{\max}(t)}  \right),
 \end{eqnarray*}
where we use that $\displaystyle\int_M e^{v(t)}dV= \int_M e^{v_0}dV$, for all $t\in [0,T)$. Integrating this last inequality, between $0$ and $t$, we get
 \begin{equation*}
 e^{v_{\max}\left(t\right)} +  \left\|Q\right\|_{L^\infty (M)} \frac{\int_{M}e^{v_0}dV}{\rho} \leq \left(e^{v_{\max}\left(0\right)}+\left\|Q\right\|_{L^\infty (M)} \frac{\int_{M}e^{v_0}dV}{\rho} \right) e^{\displaystyle \frac{\rho t}{\int_{M}e^{v_0}dV}}.
\end{equation*}
Hence \eqref{vborne1} follows.

\begin{etape}
\label{etape2}
Let $\rho\geq 8\pi$. There exists a subset $A$ of $M$, with volume $|A|>C_T$ and a constant $\delta$ depending on $M$, $Q$, $\rho$, $\left\|v_0\right\|_{C^{2+\alpha}(M)}$ and $T$ such that  
  \begin{equation}
  \label{etape2e}
  |v(x,t)|\leqslant \delta ,\ \forall x\in A\ and\ t\in [0,T).
  \end{equation}
 \end{etape}
 \bigskip
 
\noindent\textit{Proof of Step 2.} Fix $t\in [0,T)$ and set  $$M_\varepsilon = \left\{x \in M : e^{v(x,t)}<\varepsilon\right\},$$

\noindent where $\varepsilon > 0$ is a real number which will be determined later. Setting $\displaystyle\int_M e^{v_0}dV =a$, by the conservation of the volume and \eqref{vborne1}, we have :
\begin{eqnarray*}
a=\int_M e^{v(t)}dV&=&\int_{M_\varepsilon}e^{v(t)}dV + \int_{M \setminus M_\varepsilon}e^{v(t)}dV\\ 
&\leqslant &\varepsilon |M_\varepsilon | + e^{C_T} |M\setminus M_\varepsilon|.
\end{eqnarray*}
Taking $ \varepsilon = \dfrac{a}{2|M|}$, we find
\begin{equation}
\label{eg2.3}
|M \setminus M_\varepsilon | \geqslant \dfrac{a}{2}e^{-C_T}>0.
\end{equation}
Setting $A=M\backslash M_\varepsilon$, by definition of $M_\varepsilon$, we have $v(x,t)\geq \ln \varepsilon =\ln \dfrac{a}{2|M|}$, for all $x \in A$ and $t\in [0,T)$. On the other hand, by Step \ref{vborn}, $v(x,t)\leqslant C_T,\ \forall x \in M$ and $\forall t\in [0,T) ,$ therefore we get that there exists a constant $\delta$ such that
\begin{equation*}
\left|v(x,t)\right|\leqslant \delta , \hspace{12pt} \forall x \in A\ ,\ \forall t\in [0,T).
\end{equation*}

\begin{etape}
\label{etape3} 
Let $\rho\geq 8\pi$. For all $t\in [0,T)$, we have
 \begin{equation}
 \label{fardoune.g}
\int_M v^2(t)dV \leqslant C_1\int_M|\nabla v(t)|^2 dV +C_2,
\end{equation}
where $C_1$, $C_2$ are two constants depending on $T,\ Q,\ \left\|v_0\right\|_{C^{2+\alpha}(M)},\ M \ and\ A$ (where $A$ is the set defined in Step \ref{etape2}).
\end{etape}
 
\noindent\textit{Proof of Step 3.} By Poincaré's inequality, one has  
 \begin{equation}
 \label{fardoun2.87}
 \int_M v^{2}(t)dV  \leq \dfrac{1}{\lambda_1}\int_M |\nabla v(t)|^2dV+ \dfrac{1}{|M|}\left(\int_M v(t)dV\right)^2,
 \end{equation}
 where $\lambda_1$ is the first eigenvalue of the laplacian. Now, using Young's inequality and \eqref{etape2e}, we find
\begin{eqnarray}
\label{eg2.8}
\dfrac{1}{|M|}\left(\displaystyle \int_M v(t)dV \right)^2  & =& \dfrac{1}{|M|}\left(\displaystyle\int_{A}v(t)dV \right)^2 + \dfrac{1}{|M|} \displaystyle\left( \int _{M \setminus A}v (t)dV \right) ^2 \nonumber\\
 &+ &\dfrac{2}{|M|} \left( \displaystyle\int_{A}v(t)dV\right)\left(  \displaystyle\int_{M \setminus A} v(t)dV\right)\nonumber\\
&\leq &\dfrac{\delta^2 |A|^2}{|M|}+\dfrac{1}{|M|} \left(\int_{M\backslash A}v(t)dV\right)^{2} \nonumber\\
 & +& \dfrac{2\delta^2 |A|^2}{\varepsilon|M|}+\dfrac{2\varepsilon}{|M|} \left(\int_{M\backslash A}v(t)dV\right)^{2},
\end{eqnarray}
where $\varepsilon$  is a positive constant which will be determined later. By Cauchy-Schwarz's inequality, one has
 \begin{equation}
 \label{eg2.12}
 \left(\int_{M\backslash A}v(t)dV\right)^{2} \leq  |M \backslash A| \int_{M\backslash A}v^2(t)dV.
 \end{equation}
Thus, (\ref{fardoun2.87}), (\ref{eg2.8}) and (\ref{eg2.12}) yield to

 \begin{eqnarray*}
 \label{eg2.14}
 \int_M v^2(t) dV  &\leq &\frac{1}{\lambda_1}\int_M |\nabla v(t)|^2 dV+\left(1- \frac{|A|}{|M|}+\frac{2\varepsilon}{|M|}|M\backslash A|\right)\int_{M}v^2(t) dV \nonumber\\
 &+& \frac{\delta^2|A|^2}{|M|}+ \frac{2\delta^2 |A|^2}{\varepsilon|M|}.
 \end{eqnarray*}
 Choosing $\varepsilon$ such that $\alpha= \left(1- \frac{|A|}{|M|}+\frac{2\varepsilon}{|M|}|M\backslash A|\right)<1,$ and setting $\tilde{C}= \frac{\delta^2|A|^2}{|M|}+ \frac{2\delta^2 |A|^2}{|M|\varepsilon}$, we deduce that
 $$\left(1-\alpha \right)\int_M v^2 (t)dV \leq \frac{1}{\lambda_1}\int_M |\nabla v(t)|^2 dV +\tilde{C}.$$
Therefore inequality (\ref{fardoune.g}) is established.
 
\bigskip
\bigskip
\vspace{12pt}
\vspace{12pt}

\noindent \textit{Proof of Proposition \ref{etape4}.} We will consider two cases $ \rho <8\pi$ and $\rho\geq 8\pi$. In the first one, we will prove that the constant $C$ of estimate (\ref{etape4e1}) isn't depending on $T$. 
We begin with the case $\rho<8\pi$. Using  Poincaré and Young's inequalities, we have
\begin{eqnarray*}
C\int_M \left|v(t)-\bar{v}(t)\right|dV &\leq & \varepsilon \int_M |\nabla v(t)|^2dV + C,
\end{eqnarray*}
where $\varepsilon>0$ is a small constant to be chosen later. This implies that 
\begin{eqnarray}
\label{7dec2012e2}
J_\rho(v(t)) &= &\frac{1}{2}\int_{M}\left|\nabla v(t)\right|^{2}dV+\int_{M}Q \left(v(t)-\bar{v}(t)\right) dV\nonumber \\
 &-& \rho \log\left(\int_{M}e^{v(t)-\bar{v}(t)}dV\right)\nonumber\\
 &\geq &\left(\frac{1}{2}-\varepsilon \right)\int_{M}\left|\nabla v(t)\right|^{2}dV-C\nonumber\\
 &-& \rho \log\left(\int_{M}e^{v(t)-\bar{v}(t)}dV\right).
\end{eqnarray}
By Jensen's inequality, we have
\begin{equation}
\label{jensen}
\log \left(\int_M e^{v(t)-\bar{v}(t)}dV\right)\geq C,\ \forall t\in [0,T),
\end{equation}
where $\bar{v}(t)=\dfrac{\int_M v(t) dV}{|M|}$. Hence, using \eqref{7dec2012e2}, \eqref{jensen} and setting $\rho_1=\max \left(\rho ,0\right) $, we deduce that
\begin{eqnarray*}
J_\rho(v(t)) &\geq &\left(\frac{1}{2}-\varepsilon \right)\int_{M}\left|\nabla v(t)\right|^{2}dV-C\nonumber\\
 &-& \rho_1\log\left(\int_{M}e^{v(t)-\bar{v}(t)}dV\right).
\end{eqnarray*}
By Moser-Trudinger's inequality (see \cite{MR0301504}, \cite{MR0216286}), one has
\begin{equation}
\label{mosermod}
\log \int_{M}e^{\left(u(t)-\bar{u}(t)\right)}dV\leq \frac{1}{16\pi}\int_{M}\left|\nabla u(t)\right|^{2}dV+C.
\end{equation}
Therefore
\begin{eqnarray*}
\label{fardounie2.2'}
& \geq &\left(\frac{1}{2}-\frac{\rho_1}{16\pi}-\varepsilon\right)\int_{M}\left|\nabla v(t)\right|^{2}dV-C.
\end{eqnarray*}
Thus, by taking $\varepsilon= \dfrac{8\pi -\rho_1}{32\pi}$ and using that $J_\rho (v(t)) \leq J_\rho (v_0)$, $\forall t\in [0,T)$, we find that, $\forall  \rho < 8\pi $,
\begin{equation}
\label{bornegradient}
\int_{M}\left|\nabla v(t)\right|^{2}dV\leq C.
\end{equation} 
Now, using \eqref{bornegradient} and Poincaré's inequality, we obtain, $\forall \rho < 8\pi$,
\begin{equation}
\label{38bis}
 \left\|v(t) -\bar{v}(t) \right\|_{H^{1}(M)}\leq C.
 \end{equation}
\noindent Since $\displaystyle\int_{M}e^{v(t)}dV=\int_M e^{v_0}$ for all $t\in [0,T)$, using Jensen's inequality \eqref{jensen}, Moser-Trudinger's inequality \eqref{mosermod} and \eqref{bornegradient}, we deduce that
$$\left|\bar{v}(t)\right|\leq C.$$
 Finally, from (\ref{38bis}) and the previous inequality, we obtain, for all $\rho <8\pi$,
\begin{equation}
\label{ibh1}
\left\|v(t)\right\|_{H^{1}(M)}\leq C,\ \forall \ 0\leq t < T.
\end{equation}
 
Let's now consider the last case $\rho\geq 8\pi$. Since $\displaystyle\int_{M}e^{v(t)}dV=\int_M e^{v_0}$ for all $t\in [0,T)$, and by Young's inequality, we have
\begin{eqnarray}
\label{eg2.20}
J_\rho (v_0) \geq J_\rho (v)&\geq &\frac{1}{2}\int_M |\nabla v|^2 dV + \int_M Q vdV - C\nonumber \\
&\geq &\frac{1}{2}\int_M |\nabla v|^2 dV -\varepsilon \int_M v^2(t) dV-C,
\end{eqnarray}
where $\varepsilon$ is a positive constant which will be chosen later. Thanks to estimate (\ref{fardoune.g}) of Step \ref{etape3}, inequality \eqref{eg2.20} leads to
$$ \frac{1}{2}\int_M |\nabla v(t)|^2 dV \leq C + C_1 \varepsilon \int_M |\nabla v(t)|^2 dV.$$
Choosing $\varepsilon$ such that $\frac{1}{2}- \varepsilon C_1>0$, we obtain that, for all $t\in [0,T)$,
\begin{equation}
\label{fardoun2emeversionetoile}
\int_M |\nabla v(t)|^2 dV \leq C_T.
\end{equation}
Combining (\ref{fardoune.g}) and (\ref{fardoun2emeversionetoile}), we have $\displaystyle\int_M v^2(t) dV\leq C_T$. Finally, we conclude that, for all $\rho \in \R$, there exists a constant $C_T>0$ such that
$$\left\|v(t)\right\|_{H^1(M)}\leq C_T,\ \forall \ t\in [0,T).$$
\end{proof}

\begin{prop}
\label{etape5}
There exists a constant $C_T>0$ such that
$$\left\|v(t)\right\|_{H^{2}(M)}\leq C_T, \forall \ 0\leq t < T.$$
\end{prop}
\begin{proof}
Since $\left\|v\right\|_{H^1(M)}\leq C_T$, we just need to bound $\int_M (\Delta v(t))^2dV$, for all $t\in [0,T)$. To this purpose, we set
$$w(t)=\dfrac{\partial v(t)}{\partial t}e^{\frac{v(t)}{2}}.$$
By differentiating with respect to $t$ and integrating by parts on $M$, we have
\begin{eqnarray*}
&& \frac{1}{2}\dfrac{\partial }{\partial t}\int_{M}\left(\Delta v(t)\right)^{2}dV\\
& =& \int_{M}\left(w(t)e^{\frac{v(t)}{2}}+Q - \frac{\rho  e^{v(t)}}{\int_{M} e^{v_0}dV}\right) \Delta \left(w(t)e^{-\frac{v(t)}{2}}\right)dV \\
& =& -\int_{M}\left|\nabla w(t)\right|^{2}dV+\frac{1}{4}\int_{M}w^{2}(t)\left|\nabla v(t)\right|^{2}dV\\
&+&\int_{M}\Delta Q  \left(w(t)e^{-\frac{v(t)}{2}}\right)dV \\
& +&\frac{\rho}{\int_{M}e^{v_0}dV}\left(\int_{M}\nabla v(t) \nabla w(t) e^{\frac{v(t)}{2}}dV-\frac{1}{2}\int_{M}w(t) e^{\frac{v(t)}{2}}\left|\nabla v(t)\right|^{2} dV\right).
\end{eqnarray*}

\noindent Since $Q\in C^\infty (M)$ and $w(t)=\dfrac{\partial v(t)}{\partial t}e^{\frac{v(t)}{2}}$, we find
\begin{eqnarray}
\label{ve47}
&&\frac{1}{2}\frac{\partial }{\partial t}\int_{M}\left(\Delta v(t)\right)^{2}dV\nonumber\\
 & \leq &-\int_{M}\left|\nabla w(t)\right|^{2}dV+\frac{1}{4} \int_{M}w^{2}(t)\left|\nabla v(t)\right|^{2}dV+C \left\|\dfrac{\partial v(t)}{\partial t}\right\|_{L^{1}(M)}\nonumber\\
&+& C \left(\int_{M}e^{\frac{v(t)}{2}}\left(\left|\nabla w(t)\right|\left| \nabla v(t)\right|+\left|\nabla v(t)\right|^2\left|w(t)\right|\right)dV\right).
\end{eqnarray}
We now estimate the positive terms on the right of (\ref{ve47}). From Gagliardo-Nirenberg's inequality (see for example \cite{MR1961176}), $\forall f\in H^{1}(M),$
 $$\left\|f\right\|^2_{L^4(M)} \leq C\left\|f\right\|_{L^2(M)}\left\|f\right\|_{H^1(M)}, $$
using Cauchy-Schwarz's inequality and \eqref{etape4e}, we have
\begin{eqnarray}
\label{ve48}
&&\int_{M} w^2(t)\left| \nabla v(t)\right|^{2}dV\nonumber \\
&\leq & \left\|w(t)\right\|^{2}_{L^{4}(M)} \left\|\nabla v(t)\right\|^{2}_{L^{4}(M)}\nonumber \\
&\leq &C_T \left\|w(t)\right\|_{L^{2}(M)} \left\|w(t)\right\|_{H^{1}(M)} \left\| v(t)\right\|_{H^2(M)}.
\end{eqnarray}
Using \eqref{etape4e} and Moser-Trudinger's inequality \eqref{mosermod}, we deduce that there exists a constant $C_T$ such that, for all $ t\in [0,T)$,  and $p\in \R$,
\begin{equation}
\label{ibep}
\int_M e^{pv(t)}dV\leq C_T.
\end{equation}
In the same way as for proving \eqref{ve48}, using \eqref{etape4e} and \eqref{ibep}, we have
\begin{eqnarray}
\label{ve50}
&&\displaystyle \int_{M} \left|\nabla v(t)\right|^{2} \left| w(t)\right|e^{\frac{v(t)}{2}}dV\nonumber \\
 & \leq &\left(\displaystyle\int_{M}\left|\nabla v(t)\right|^{4}dV\right)^{\frac{1}{2}} \left(\displaystyle\int_{M}w^{4}(t)dV\right)^{\frac{1}{4}}  \left(\displaystyle\int_{M}e^{2v(t)}dV\right)^{\frac{1}{4}}\nonumber \\
& \leq & C_T \left\|v(t)\right\|_{H^2(M)}\left\|w(t)\right\|_{H^{1}(M)}^{\frac{1}{2}}\left\|w(t)\right\|_{L^{2}(M)}^{\frac{1}{2}},
\end{eqnarray}
\begin{eqnarray}
\label{ve51}
&&\int_{M} \left|\nabla w(t)\right| \left| \nabla v(t)\right|e^{\frac{v(t)}{2}}dV \nonumber \\
& \leq &\left(\int_{M}\left|\nabla w(t)\right|^{2}dV\right)^{\frac{1}{2}} \left(\int_{M}\left|\nabla v(t)\right|^{4}dV\right)^{\frac{1}{4}}  \left(\int_{M}e^{2v(t)}dV\right)^{\frac{1}{4}}\nonumber \\
& \leq &C_T \left\|w(t)\right\|_{H^{1}(M)}\left\|v(t)\right\|_{H^2(M)}^{\frac{1}{2}},
\end{eqnarray}
and
\begin{eqnarray}
\label{51ter2}
\int_M \left|\dfrac{\partial v(t)}{\partial t}\right|dV &\leq &\left(\int_M \left(\dfrac{\partial v(t)}{\partial t}\right)^{2}e^{v(t)}dV\right)^{\frac{1}{2}}\left(\int_M e^{-v(t)}dV\right)^{\frac{1}{2}}\nonumber \\
&\leq & C_T \left\|w(t)\right\|_{L^2(M)}.
\end{eqnarray}
Finally, putting (\ref{ve48}), (\ref{ve50}), (\ref{ve51}), (\ref{51ter2}) in (\ref{ve47}), we obtain
\begin{eqnarray*}
&&\frac{1}{2}\dfrac{\partial }{\partial t}\displaystyle\int_{M}\left(\Delta v(t)\right)^{2}dV \\
& \leq &- \displaystyle\int_M |\nabla w(t)|^2dV+C_T \left\|w(t)\right\|_{H^{1}(M)}\left\|w(t)\right\|_{L^{2}(M)}\left\|v(t)\right\|_{H^2(M)}+C_T\left\|w(t)\right\|_{L^{2}(M)} \\
& +& C_T \left(\left\|w(t)\right\|_{H^{1}(M)}\left\|v (t)\right\|_{H^2(M)}^{\frac{1}{2}}+ \left\|w(t)\right\|_{H^{1}(M)}^{\frac{1}{2}}\left\|w(t)\right\|_{L^{2}(M)}^{\frac{1}{2}}\left\|v(t)\right\|_{H^2(M)}\right).
\end{eqnarray*}
Using Young's inequality, we get
\begin{eqnarray}
\label{18juin2012}
&&\dfrac{\partial }{\partial t}\left(\int_{M}\left(\Delta v(t)\right)^{2}dV+ 1\right)  \nonumber \\
&\leq & C_T \left(\int_{M}\left(\Delta v(t)\right)^{2}dV+ 1\right)\left(\left\|w(t)\right\|_{L^{2}(M)}^{2}+1\right).
\end{eqnarray}
\noindent On the other hand, by (\ref{eg2.5'}), we have, for all $t\in [0,T)$,
\begin{eqnarray}
\label{51bis}
\int_0^t\left\|w(s)\right\|_{L^{2}(M)}^{2}ds &=&\int_0^t\int_{M}\left(\dfrac{\partial v(s)}{\partial s}\right)^{2}e^{v(s)}dVds\nonumber \\
&=&-\int_0^t\dfrac{\partial}{\partial s} J_\rho(v(s))ds=J_\rho(v_0)-J_\rho(v(t))\nonumber \\
 &\leq &C_T,
\end{eqnarray}
where we use the fact that $\left\|v(t)\right\|_{H^1(M)}\leq C_T$ from Proposition \ref{etape4}. Integrating \eqref{18juin2012} with respect to $t$ and using (\ref{51bis}), we have
\begin{eqnarray*}
\int_M \left(\Delta v(t)\right)^2 dV &\leq & C_T,\ \forall t\in[0,T) .
\end{eqnarray*}
Since $\left\|v(t)\right\|_{H^1(M)}\leq C_T$, we deduce that
$$\left\|v(t)\right\|_{H^2(M)}\leq C_T, \forall \ t\in \left[0,T\right).$$
\end{proof}

\bigskip
\bigskip

\noindent\textit{Proof of Theorem \ref{exglob}.} We recall that to prove the global existence of the flow it is sufficient to prove \eqref{7dec2012e1}, i.e. there exists a constant $C_T$ depending on $T$ such that
$$\left\|v\right\|_{C^{2+\alpha,1+\frac{\alpha}{2}}(M\times [0,T))}\leq C_T,\ \alpha\in (0,1).$$   
First, we claim that for all $\alpha \in (0,1)$, there exists a constant $C_T$ such that
\begin{equation}
\label{brenreg}
|v(x_1,t_1)-v(x_2,t_2)|\leq C_T(|t_1-t_2|^{\frac{\alpha}{2}}+|x_1-x_2|^{\alpha}),
\end{equation}
for all $x_1,x_2 \in M$, all $t_1,t_2 \in [0,T)$, and where $|x_1-x_2|$ stands for the geodesic distance from $x_1$ to $x_2$ with respect to the metric $g$. From Proposition \ref{etape5}, we have, for all $t\in [0,T)$,  $\left\|v(t)\right\|_{H^2(M)}\leq C_T$. Thus, by Sobolev's embedding Theorem (see \cite{hebey}), we get that, for $\alpha \in (0,1)$, $v(t)\in C^\alpha (M)$ that is, for all $x,y\in M$,
\begin{equation}
\label{bren4}
|v(x,t)-v(y,t)|\leq C_T |x-y|^\alpha .
\end{equation}
If $t_2-t_1 \geq 1$, using \eqref{bren4}, it is easy to see that \eqref{brenreg} holds. Thus, from now on, we assume that $0<t_2-t_1<1$. Since $v(t)$ is a solution of \eqref{E:flot} and $\left\|e^{v(t)}\right\|_{C^\alpha (M)}\leq C_T$, for all $t\in [0,T)$, one has
$$\left|\dfrac{\partial v(t)}{\partial t}\right|^2\leq C_T |\Delta v(t)|^2 +C_T.$$
Integrating on $M$, we obtain, for all $t \in [0,T)$,
\begin{equation}
\label{bren1}
\int_M \left|\dfrac{\partial v(t)}{\partial t}\right|^2dV \leq C_T \left\|v(t)\right\|^2_{H^2(M)}+C_T\leq C_T.
\end{equation}
Now, we write 
\begin{eqnarray}
\label{bren3}
|v(x,t_1)&-&v(x,t_2)|= \frac{1}{|B_{\sqrt{t_2-t_1}}(x)|}\int_{B_{\sqrt{t_2-t_1}}(x)} |v(x,t_1)-v(x,t_2)|dV(y)\nonumber \\
&\leq & \frac{C}{ t_2-t_1} \int_{B_{\sqrt{t_2-t_1}}(x)} |v(x,t_1)-v(y,t_1)|dV(y)\nonumber\\
& +&\frac{C}{ t_2-t_1}\int_{B_{\sqrt{t_2-t_1}}(x)} |v(y,t_1)-v(y,t_2)|dV(y) \nonumber\\
&+&\frac{C}{ t_2-t_1}\int_{B_{\sqrt{t_2-t_1}}(x)} |v(y,t_2)-v(x,t_2)|dV(y),
\end{eqnarray}
where $B_{\sqrt{t_2-t_1}}(x)$ stands for the geodesic ball of center $x$ and radius $\sqrt{t_2-t_1}$. Let's consider the first term on the right of (\ref{bren3}). Using (\ref{bren4}), we obtain
\begin{eqnarray}
\label{bren5}
&&\frac{C}{t_2-t_1} \int_{B_{\sqrt{t_2-t_1}}(x)} |v(x,t_1)-v(y,t_1)|dV(y)\nonumber \\ &\leq & \frac{C_T}{ (t_2-t_1)} \int_{B_{\sqrt{t_2-t_1}}(x)} |x-y|^{\alpha}dV(y)\nonumber\\
&\leq & C_T(t_2-t_1)^{\frac{\alpha}{2}}.
\end{eqnarray}
In the same way, we have
\begin{equation}
\label{bren6}
\frac{C}{ t_2-t_1} \int_{B_{\sqrt{t_2-t_1}}(x)} |v(x,t_2)-v(y,t_2)|dV(y) \leq C_T (t_2-t_1)^{\frac{\alpha}{2}}.
\end{equation} 
We find, using H\"{o}lder's inequality and (\ref{bren1}),
\begin{eqnarray}
\label{bren2}
&&\frac{C}{t_2-t_1}\int_{B_{\sqrt{t_2-t_1}}(x)} |v(y,t_1)-v(y,t_2)|dV(y)\nonumber \\ &\leq & C \sup_{t_1\leq\tau \leq t_2 }\int_{B_{\sqrt{t_2-t_1}}(x)} \left|\dfrac{\partial v}{\partial s}  \right| (y,\tau)dV(y)\nonumber \\
& \leq & C \sqrt{t_2-t_1} \sup_{t_1\leq\tau \leq t_2 }\left(\int_{B_{\sqrt{t_2-t_1}}(x)} \left|\dfrac{\partial v}{\partial s}  \right|^2 (y,\tau)dV(y)\right)^{\frac{1}{2}}\nonumber \\
& \leq & C_T \sqrt{t_2-t_1}.
\end{eqnarray}
Putting (\ref{bren5}), (\ref{bren6}), (\ref{bren2}) in (\ref{bren3}) and noticing that for all $0<t_2-t_1<1$, we have $\sqrt{t_2-t_1}\leq (t_2-t_1)^{\frac{\alpha}{2}}$, we find
\begin{equation}
\label{bren7}
|v(x,t_1)-v(x,t_2)|\leq  C_T (t_2-t_1)^{\frac{\alpha}{2}} .
\end{equation}
Therefore, from (\ref{bren4}) and (\ref{bren7}), we see that \eqref{brenreg} holds. 
In view of \eqref{brenreg}, we may apply the standard regularity theory for parabolic equations (see for example \cite{MR0181836}) to derive that there exists a constant $C_T$ depending on $T$ such that
$$\left\|v\right\|_{C^{2+\alpha,1+\frac{\alpha}{2}}(M\times [0,T))}\leq C_T,\ \alpha\in (0,1).$$
This completes the existence part of Theorem \ref{exglob}. The uniqueness follows from Proposition \ref{propdonin}. 
\begin{rmq}
\label{exglob26mai}
Following the proof of Theorem \ref{exglob}, we see that, for all $T>0$ fixed, if $\left\|u_0\right\|_{C^{2+\alpha}(M)}\leq K$ for some constant $K>0$, then there exists a constant $C_T>0$ depending on $K$ and $T$ such that
$$\left\|u\right\|_{C^{2+\alpha,1+\frac{\alpha}{2}}(M\times [0,T])}\leq C_T .$$
\end{rmq}
\subsection{Continuity of the flow with respect to its initial data.}
We now show the continuity of the flow with respect to its initial data which will be useful for the proof of Theorem \ref{conv2} (see Section $3$).
\begin{prop}
\label{propdonin}
Let $u,v \in C^{2+\alpha,1+\frac{\alpha}{2}}_{loc}(M\times [0,+\infty ))$, $\alpha\in (0,1)$, be solutions of 
\begin{equation*}
\left\{\begin{array}{ll}
 \dfrac{ \partial}{\partial t}e^{v} =\bigtriangleup v - Q+ \rho \dfrac{  e^{v}}{\int_M  e^{v}dV} \\
 v(x,0)=v_{0}(x),
\end{array}
\right.
\end{equation*}
and
\begin{equation*}
\left\{\begin{array}{ll}
 \dfrac{ \partial}{\partial t}e^{u} =\bigtriangleup u - Q+ \rho \dfrac{  e^{u}}{\int_M  e^{u}dV} \\
 u(x,0)=u_{0}(x),
\end{array}
\right.
\end{equation*}
where $u_0 ,v_0 \in C^{2+\alpha}(M)$. Then, for all $T>0$, there exists a constant $C_T>0$ depending on $\left\|u_0\right\|_{C^{2+\alpha}(M)},\left\|v_0\right\|_{C^{2+\alpha}(M)}$ and $T$, such that
\begin{equation}
\label{10nov2012e1}
\left\|u -v\right\|_{C^{2+\alpha,1+\frac{\alpha}{2}}(M\times [0,T])}\leq C_T \left\|u_0-v_0\right\|_{C^{2+\alpha}(M)} .
\end{equation}
\end{prop}
\begin{rmq}
\label{propdonin26mai}
Following the proof of Proposition \ref{propdonin}, we see that, for all $T>0$ fixed, if $\left\|u_0\right\|_{C^{2+\alpha}(M)}\leq K_1$ and $\left\|v_0\right\|_{C^{2+\alpha}(M)}\leq K_2$ for some constants $K_1,K_2>0$, then there exists a constant $C_T>0$ depending on $K_1$, $K_2$ and $T$ such that
$$\left\|u -v\right\|_{C^{2+\alpha,1+\frac{\alpha}{2}}(M\times [0,T])}\leq C_T \left\|u_0-v_0\right\|_{C^{2+\alpha}(M)} .$$
\end{rmq}
\begin{proof}
We set $F=u-v$. F satisfies
$$\frac{\partial F(x,t)}{\partial t}=a(x,t)\Delta F(x,t)+b(x,t)F(x,t)+f,$$
where $$a(x,t)=\frac{1}{2}(e^{-u(x,t)}+e^{-v(x,t)}),$$ $$b(x,t)=-\left(\frac{1}{2}\left(\Delta u(x,t)+\Delta v(x,t)\right)-Q(x)\right)\int_0^1 e^{-su(x,t)-(1-s)v(x,t)}ds$$ and $$f=\rho \left(\frac{1}{\int_M e^{u_0}dV}-\frac{1}{\int_M e^{v_0}dV}\right).$$ Since $u,v \in C^{2+\alpha,1+\frac{\alpha}{2}}(M\times [0,T])$, there exists two constants $C_T$ such that 
\begin{equation}
\label{propdoninest1}
\left\|a\right\|_{C^{\alpha,\frac{\alpha}{2}}(M\times [0,T])}\leq C_T,
\end{equation}
and
\begin{equation}
\label{propdoninest2}
\left\|b\right\|_{C^{\alpha,\frac{\alpha}{2}}(M\times [0,T])}\leq C_T.
\end{equation}
It is easy to check that
\begin{equation}
\label{donin1}
|f|\leq C \left\|F(.,0)\right\|_{C^{2+\alpha}(M)},
\end{equation}
where $C$ is a constant depending on $\left\|u_0\right\|_{C^{2+\alpha}(M)}$ and $\left\|v_0\right\|_{C^{2+\alpha}(M)}$. Using (\ref{propdoninest1}), (\ref{propdoninest2}) and regularity theory for parabolic equations (see \cite{MR0181836}), we find that, for all $T>0$,
\begin{equation}
\label{donin2}
\left\|F\right\|_{C^{2+\alpha,1+\frac{\alpha}{2}}(M\times [0,T])}\leq C_T \left( \left\|F(.,0)\right\|_{C^{2+\alpha}(M)}+\left\|F\right\|_{L^\infty (M\times [0,T])}\right).
\end{equation}
Let's set $F_{\max}(t)=\displaystyle\max_{x\in M} F(x,t)$ and $F_{\min} (t)=\displaystyle\min_{x\in M} F(x,t)$ for $t\in [0,T]$. By the maximum principle and using (\ref{donin1}), we have
$$F_{\max}(t)\leq e^{T \left\|b\right\|_{L^\infty (M\times[0,T])}} \left|F_{\max}(0)\right|+ C_T \left\|F(.,0)\right\|_{C^{2+\alpha}(M)}e^{T \left\|b\right\|_{L^\infty (M,[0,T])}}$$
and
$$F_{\min}(t)\geq -e^{T \left\|b\right\|_{L^\infty (M\times[0,T])}} \left|F_{\min}(0)\right|- C_T \left\|F(.,0)\right\|_{C^{2+\alpha}(M)}e^{T \left\|b\right\|_{L^\infty (M,[0,T])}}.$$
This yields to
\begin{equation}
\label{donin3}
\left\|F\right\|_{L^\infty (M\times [0,T])}\leq C_T\left\|F(.,0)\right\|_{C^{2+\alpha}(M)}.
\end{equation} 
Combining (\ref{donin2}) and (\ref{donin3}), we obtain \eqref{10nov2012e1},
$$\left\|u -v\right\|_{C^{2+\alpha,1+\frac{\alpha}{2}}(M\times [0,T])}\leq C_T \left\|u_0-v_0\right\|_{C^{2+\alpha}(M)} .$$

\end{proof}
\section{Compactness property.}
In this section, we study the compactness property of $(v_n)_n \subseteq H^2(M)$ solutions of \eqref{regjanveqini}
\begin{equation*}
-\Delta v_n =Q_0 +h_n e^{v_n}+\rho e^{v_n},
\end{equation*}
where $\rho>0$, $Q_0 \in C^0 (M)$ is a given function and $(h_n)_n \subseteq C^0(M)$.
\subsection{Concentration-compactness.}
We first give a concentration-compactness criterion in the spirit of Brezis-Merle \cite{MR1132783} for solutions $(v_{n})_n\subseteq H^2(M)$ of 
$$-\Delta v_{n}=F_{n},$$
where $(F_{n})_{n}\subseteq L^{1}(M)$ is a sequence of functions such that $\displaystyle\int_{M}|F_{n}|dV\leq C$ for all $n\geq 0$. 
\begin{prop}
\label{mp1}
 Let $(v_n)_n \subseteq H^2(M)$ be a sequence of solutions of $-\Delta v_{n}=F_{n}$ with $(F_{n})_{n}\subseteq L^{1}(M)$ such that $\displaystyle\int_{M}|F_{n}|dV\leq C$ for all $n\geq 0$. Then, we have up to a subsequence :  \newline
-(i) either, there exist constants $C>0$ and $\alpha>1$ such that
 $$\int_{M}e^{\alpha (v_{n}-\bar{v}_{n})}dV\leq C,\ \forall \ n\geq 0.$$ 
 -(ii) or, there exist points $x^{1},...,x^{L}\in M$ such that, $\forall r>0$ and $\forall i\in \left\{1,...,L\right\}$, we have 
 $$\liminf_{n \rightarrow \infty} \int_{B_{r}(x_{i})} |F_{n}|dV\geq 4\pi.$$  
 \end{prop}
\begin{proof}
The proof is well known. We give it for the sake of selfcontainedness. Let's suppose that (ii) doesn't hold, i.e. 
 \begin{equation}
 \label{me1}
 \forall x \in M, \exists \ r_{x}>0 \ such\ that\ \int_{B_{r_{x}}(x)}|F_{n}|dV\leq 4\pi -\delta_{x},
\end{equation} 
where $\delta_{x}>0$ and $n$ large enough. Since M is compact, we can cover it by a finite number $j$ of balls $B_{i}:=B_{\frac{r_{x^{i}}}{2}}(x^{i}),\ i\in \left\{1,...,j\right\}$. We set $B_{r_{x^{i}}}(x^{i})=\tilde{B_{i}}$. Thus we have, for all $x\in B_{i}$, 
 $$v_{n}(x)-\bar{v}_{n}=\int_{\tilde{B_{i}}}F_{n}(y) G(x,y) dV(y)+ \int_{M\backslash \tilde{B_{i}}} F_{n}(y) G(x,y) dV(y),$$
 where $G:M\times M\rightarrow \R$ is the Green function of the laplacian (see for example \cite{MR1636569}). Multiplying by $\alpha>0$ and taking the exponential, we find,  $\forall x \in B_{i}$,
 \begin{eqnarray}
 \label{me2}
 \exp\left[\alpha (v_{n}(x)-\bar{v}_{n})\right]&=&\exp\left[\int_{\tilde{B_{i}}}\alpha F_{n}(y) G(x,y) dV(y)\right]\nonumber \\
  &\times & \exp\left[ \int_{M\backslash \tilde{B_{i}}}\alpha F_{n}(y) G(x,y) dV(y)\right].
 \end{eqnarray}
Since $G$ is smooth outside its diagonal and $\int_{M} |F_{n}|dV\leq C$, there exists a constant $C$ such that, for all $x\in B_i$,
 $$\exp\left[ \int_{M\backslash \tilde{B_{i}}}\alpha F_{n}(y) G(x,y) dV(y)\right]\leq C.$$
 Combining the previous inequality, (\ref{me2}), and integrating on $B_{i}$, we get
\begin{eqnarray}
\label{mp1.1}
&&\int_{B_{i}}\exp\left[\alpha (v_{n}(x)-\bar{v}_{n})\right]dV(x)\nonumber \\
&\leq &C \int_{B_{i}}\exp\left[\int_{\tilde{B_{i}}}\alpha F_{n}(y) G(x,y) dV(y)\right]dV(x).
\end{eqnarray}
Using Jensen's inequality, we have
\begin{eqnarray*}&&\exp\left[\int_{\tilde{B_{i}}}\alpha |F_{n}(y)| |G(x,y)| dV(y)\right] \\
&\leq &\int_{\tilde{B_{i}}} \exp\left[ \alpha \left\|F_{n} \right\|_{L^{1}(\tilde{B_{i}})} |G(x,y)|\right]\frac{|F_{n}|(y)}{\left\|f_{n} \right\|_{L^{1}(\tilde{B_{i}})}}dV(y).
\end{eqnarray*}
Thus, using Fubini's Theorem in (\ref{mp1.1}), we obtain 
$$\int_{B_{i}}\exp\left[\alpha (v_{n}(x)-\bar{v}_{n})\right]dV(x)\leq C \sup_{y\in \in \tilde{B}_i } \int_{B_i} \exp\left[ \alpha \left\|F_{n} \right\|_{L^{1}(\tilde{B_{i}})} |G(x,y)|\right]dV(x).$$
Since $\left|G(x,y) -\frac{1}{2\pi}\log \frac{1}{|x-y|}\right|\leq C$, for all $x\neq y\in M$, and since $\left\|F_{n} \right\|_{L^{1}(\tilde{B_{i}})}\leq C$, we finally arrive at
$$\int_{B_{i}}\exp\left[\alpha (v_{n}(x)-\bar{v}_{n})\right]dV(x)\leq C \int_{M}\left(\frac{1}{|x-y|}\right)^{\frac{\alpha \left\|F_{n} \right\|_{L^{1}(\tilde{B_{i}})}}{2\pi}}dV(x).$$
Using (\ref{me1}), we see that there exists $\alpha>1$ such that
$$ \alpha \int_{\tilde{B_{i}}}|F_{n}|dV<4\pi.$$
Thus, we have
$$\int_{B_{i}}e^{\alpha (v_{n}-\bar{v}_{n})}dV< +\infty \ ,\ \forall i\in \{1,...,j\}.$$
Since $M$ is covered by a finite number of $B_i$, we obtain 
$$\int_{M}e^{\alpha (v_{n}-\bar{v}_{n})}dV\leq C.$$
\end{proof} 
Now, we prove that if the first alternative of the previous proposition occurs then the solutions of (\ref{regjanveqini}) are compact. 

\begin{prop}
\label{pasblow}
Let $(v_n)_n \subseteq H^2 (M)$ be a sequence of solutions of (\ref{regjanveqini}) such that $\displaystyle\int_M e^{v_n}dV=1$, $\forall n\geq 0$ and satisfying conditions \eqref{regdeceq1}. Suppose there exists constants $C>0$ and $\alpha>1$ such that
$$\int_M e^{\alpha (v_n-\bar{v}_n)}dV\leq C.$$
Then there exists a constant $C$ such that
$$\left\|v_n\right\|_{H^2(M)}\leq C,\ \forall n\geq 0.$$
\end{prop}

\begin{proof}
From $\displaystyle\int_M e^{v_n}dV=1$ for all $n\geq 0$ and (\ref{regdeceq1})(i), we have $\displaystyle\lim_{n\rightarrow +\infty}\int_M \left|h_n\right|\ e^{v_n}dV=0.$ Thus, there exists a constant $C$ such that
$$\left\|\rho e^{v_n}+Q_0+h_n e^{v_n}\right\|_{L^{1}(M)}\leq C.$$
Therefore we can apply Proposition \ref{mp1} for $F_n=\rho e^{v_n}+Q_0+h_n e^{v_n}$. Since $\displaystyle\int_M  e^{\alpha (v_n-\bar{v}_n)}dV\leq C$, using H\"{o}lder's inequality, we obtain
$$1=e^{\bar{v}_n}\int_M e^{(v_n-\bar{v}_n)}dV \leq C e^{\bar{v}_n}\left(\int_M  e^{\alpha (v_n-\bar{v}_n)}dV\right)^{\frac{1}{\alpha}}\leq Ce^{\bar{v}_n} .$$
From the previous inequality and Jensen's inequality, we deduce that there exists a constant $C$ such that
\begin{equation}
\label{pasblow1}
\int_M e^{\alpha v_n}dV \leq C.
\end{equation}
Using H\"{o}lder's inequality and \eqref{pasblow1}, we have, for $\gamma = \frac{2\alpha}{\alpha +1}>1$,
\begin{eqnarray}
\label{pasblow2}
\int_M |h_n e^{v_n}|^\gamma dV & \leq &\left(\int_M h_n^2 e^{v_n}dV\right)^{\frac{\gamma}{2}} \left(\int_M e^{\alpha v_n}dV\right)^{1-\frac{\gamma}{2}}\nonumber \\
&\leq & C.
\end{eqnarray}
We also have, since $\gamma <\alpha$,
\begin{equation}
\label{pasblow4}
\int_M e^{\gamma v_n}dV\leq C\left(\int_M e^{\alpha v_n}dV\right)^{\frac{\gamma}{\alpha}}\leq C.
\end{equation}
We deduce, from (\ref{pasblow2}) and (\ref{pasblow4}), that $\left\|F_n\right\|_{L^{\gamma}(M)}\leq C.$ Since $(v_n)_n$ satisfies $-\Delta v_n=F_n$, by elliptic regularity and Sobolev's embedding Theorem, there exist two constants $C$ and $\beta \in (0,1)$ such that 
\begin{equation}
\label{pasblow3}
\left\|v_n\right\|_{C^{\beta}(M)}\leq C.
\end{equation}
Using (\ref{pasblow3}), we obtain that $\left\|e^{v_n}\right\|_{C^{\beta}(M)}\leq C$ and using estimate \eqref{regdeceq1}(i), we get
\begin{eqnarray*}
\int_M |h_n e^{v_n}|^2 dV &\leq &C \left\|e^{v_n}\right\|_{C^{\beta}(M)} \int_M h_n^2 e^{v_n}dV
\leq  C.
\end{eqnarray*}
Hence $\left\|F_n\right\|_{L^2(M)}\leq C$. Therefore, by elliptic regularity, there exists a constant $C$ such that $$\left\|v_n \right\|_{H^2(M)}\leq C.$$  
\end{proof}
In the sequel, we will say that solutions $(v_n)_n \subseteq H^2(M)$ of \eqref{regjanveqini} such that $\displaystyle\int_M e^{v_n}dV=1$, $\forall n\geq 0$, and satisfying \eqref{regdeceq1} are non-compact if the first alternative of Proposition \ref{mp1} doesn't occur.

\subsection{Asymptotic profil}
We will study the asymptotic behaviour of non-compact solutions of equation \eqref{regjanveqini}. We follow the ideas of Malchiodi \cite{MR2248155} for the study of compactness for solutions to a fourth-order equation. Thanks to Proposition \ref{mp1}, we can show that if $(v_n)_n \subseteq H^2(M)$ is a sequence of non-compact solutions of \eqref{regjanveqini}, then there exists $\beta \in \left(0, \frac{\pi}{\rho}\right)$, radii $(r_{n})_{n},(\hat{r}_{n})_{n}$ and points $(x_{n})_{n}\subseteq M$ satisfying the following properties 
\begin{eqnarray}
\label{me32}
&&\hat{r}_{n}\underset{n\rightarrow +\infty}{\longrightarrow} 0\ ;\ \frac{r_{n}}{\hat{r}_{n}}\underset{n\rightarrow +\infty}{\longrightarrow} 0\ 
;\ \int_{B_{r_{n}}(x_{n})} e^{v_{n}}dV=\beta ;\nonumber \\
&& \ \int_{B_{r_{n}}(y)} e^{v_{n}}dV
< \frac{\pi}{\rho},\ \forall y \in B_{\hat{r}_{n}}(x_{n}).
\end{eqnarray}
Now, in order to understand the asymptotic behavious of $(v_n)_n$ near $(x_n)_n$, we introduce the following tranformation of $(v_n)_n$ 
$$\hat{v}_{n}=v_n(\exp_{x_n}(r_n .))+2\log r_n, $$
where $\exp_{x_n}$ stands for the application exponential centered in $x_n$. We obtain the following proposition.
\begin{prop} 
\label{mp2}
Let $(v_n)_n \subseteq H^2(M)$ be a sequence of non-compact solutions of \eqref{regjanveqini} such that $\displaystyle\int_M e^{v_n}dV=1$, $\forall n\geq 0$ and satisfying \eqref{regdeceq1}. Then there exist a sequence of points $(x_n)_{n}$ and a sequence of real numbers $(r_n)_n$ such that
$$\hat{v}_{n}=v_n(\exp_{x_n}(r_n .))+2\log r_n \underset{n\rightarrow +\infty}{\longrightarrow} \hat{v}_{\infty}\ \ in\ C_{loc}^{\alpha}(\R^{2}),\ \alpha\in (0,1),$$
and weakly in $H^2_{loc}(\R^2)$ where $\hat{v}_{\infty}\in C^\infty (\R^2)$ is solution of the equation $$-\Delta_{g_{\R^{2}}}\hat{v}_{\infty}=\rho e^{\hat{v}_{\infty}},$$
with $\displaystyle\int_{\R^2}e^{\hat{v}_\infty}dV_{g_{\R^2}}\leq 1$. Moreover, for any sequence $(\beta_n)_n$, $\beta_{n}\underset{n\rightarrow +\infty}{\longrightarrow} +\infty$, there exists a sequence $(b_n)_n$, $b_{n}\underset{n\rightarrow +\infty}{\longrightarrow} +\infty$, such that $b_n \leq \beta_n$, $\forall n\geq 0$ and 
\begin{equation}
\label{me39}
\lim_{n\rightarrow +\infty}\int_{B_{b_{n}r_{n}}(x_{n})}e^{v_{n}}dV= \frac{8\pi }{\rho}.
\end{equation}
We will say that the sequence $(b_n)_n$ tends arbitrarily slowly to $+\infty$.
\end{prop}
\begin{rmq}
Let $\hat{v}_\infty\in C^\infty (\R^2)$ be defined as in the previous proposition. Then, the Classification Theorem of Chen-Li \cite{MR1121147} implies that there exists $\lambda>0$ and $x_{0}\in \R^{2}$  such that 
$$\hat{v}_{\infty}(x)=2\log\frac{2\lambda}{1+(\lambda |x-x_{0}|)^{2}}+\log\frac{2}{\rho}.$$
\end{rmq}
\begin{proof}
\setcounter{etape}{0}
The proof follows the same steps than Malchiodi \cite{MR2248155}. So we omit it.
\end{proof}
In view of the previous proposition, we now define the notion of blow-up sequence for $(v_n)_n \subseteq H^2(M)$ solutions of \eqref{regjanveqini}.
\begin{defi}
\label{regdecdef1}
Let $(x_n)_n$ be a sequence of points of $M$ and $(r_n)_n$ a sequence of positive real numbers such that $r_n \underset{n\rightarrow +\infty}{\longrightarrow} 0$. We say that $(x_n,r_n)_n$ is a blow-up for $(v_n)_n$ solutions of \eqref{regjanveqini} if the sequence defined by $\hat{v}_n (x)=v_n (\exp_{x_n}(r_n x))+2\log r_n$ where $x\in \R^2$, $|x|<\frac{\delta}{r_n}$ and $\delta< i(M)$ (where $i(M)$ stands for the injectivity radius of $M$), converges in $C^\alpha (B_R^{\R^2})$, $(\alpha \in (0,1))$, for all $R>0$ fixed, to a solution $\hat{v}_\infty \in C^\infty (\R^2)$ of equation 
$$-\Delta_{\R^2}\hat{v}_\infty =\rho e^{\hat{v}_\infty},$$
with $\displaystyle\int_{\R^2}e^{\hat{v}_n}dV_{\R^2}<+\infty$. 
\end{defi}
\subsection{Integral Harnack type inequality.}
Now, we will prove an integral Harnack type inequality for solutions $(v_n)_n \subseteq H^2(M)$ of equation $-\Delta v_n =Q_0+f_n$ where $(f_n)_n\subseteq L^1(M)$ and $Q_0\in C^0 (M)$. In the following, $f_n^+=\max \left(0, f_n\right)$ (resp. $f_n^-=-\min\left(0,f_n \right)$) stands for the positive (resp. negative) part of $f_n$.\newline
We obtain the following integral Harnack type inequality.
\begin{prop}
Let $(v_n)_n \subseteq H^2(M)$ and $(f_n)_n\subseteq L^1(M)$ be two sequences of functions such that 
\begin{equation}
\label{11nov2012e1}
-\Delta v_n =Q_0+f_n,
\end{equation}
with $Q_0\in C^0 (M)$, $\left\|f_n\right\|_{L^1(M)}\leq C_1$ and suppose that there exists $p>1$ such that $\left\|f_n^{-}\right\|_{L^p(M)}\leq C_2$, where $C_1$ and $C_2$ are two constants not depending on $n$. Then, for all $\beta >0$, there exists a positive constant $C_\beta$ depending on $\beta$, $C_1$ and $C_2$ such that for all $x_0,\ y_0\in M$ and $r,R \in (0,i(M))$ satisfying $r\leq R$, $|x_0-y_0|\leq \beta R$ and $\displaystyle\int_{B_{2R}(x_0)}f_n^{+}dV\leq 2\pi$, we have
$$\int_{B_R(x_0)}e^{v_n}dV \leq C_\beta \left(\dfrac{r}{R}\right)^{\frac{\left\|f_n^{+}\right\|_{L^1(B_r (y_0))}}{2\pi}-2}\int_{B_r(y_0)}e^{v_n}dV.$$ 
\label{regdecthm1}
\end{prop}
\begin{proof} 
We begin with some remarks on the Green function $G(.,.): M\times M\rightarrow \R$ of the laplacian. Up to adding a constant to $G$, we can assume without loss of generality (see for example \cite{MR1636569}) that $G(x,y)\geq 0$ for all $x,y \in M$. Note also that for all $q\in \left[1,+\infty\right)$ and $x\in M$, $\left\|G(x,.)\right\|_{L^q(M)}\leq C_q$, where $C_q$ is a constant depending on $M$ and $q$. We will divide the proof of Proposition \ref{regdecthm1} into $2$ steps. 
\setcounter{etape}{0}
\begin{etape}
\label{regdecthm1etape1}
We have
\begin{eqnarray}
\label{regdecthm1etape1.1}
\int_{B_R(x_0)}e^{v_n(x)}dV(x)
&\leq &\exp \left( C+\bar{v}_n +\displaystyle\sup_{x\in B_R(x_0)}\int_{M\backslash B_{2R}(x_0)} G(x,y) f_n^+ dV(y)\right)\nonumber \\
&\times & R^{2-\frac{\left\|f_n^{+}\right\|_{L^1(B_{2R} (x_0))}}{2\pi}}.
\end{eqnarray}
\end{etape}
\noindent\textit{Proof of Step 1.} Let $x\in B_R(x_0)$. Since $(v_n)_n \subseteq H^2(M)$ is a sequence of solutions of \eqref{11nov2012e1} and $\left\|G(x,.)\right\|_{L^1(M)}\leq C$, using Green's representation formula, we have
\begin{eqnarray*}
v_n(x)&=& \bar{v}_n +\int_M G(x,y) (f_n(y)+Q_0(y))dV(y)\\
&\leq & \bar{v}_n +\int_{B_{2R}(x_0)}G(x,y) f_n^+(y) dV(y)\\
&+&\sup_{x\in B_R(x_0)}\int_{M\backslash B_{2R}(x_0)} G(x,y) f_n^+(y) dV(y)+C. 
 \end{eqnarray*}  
Thus integrating on $B_R(x_0)$, we find
\begin{eqnarray}
\label{regdecthm1etape1.2}
\int_{B_R(x_0)}e^{v_n(x)}dV(x)&\leq & \exp \left(C+\bar{v}_n +\displaystyle\sup_{x\in B_R(x_0)}\int_{M\backslash B_{2R}(x_0)} G(x,y) f_n^+(y) dV(y)\right) \nonumber\\
&\times & \int_{B_R (x_0)}\exp \left(\int_{B_{2R}(x_0)}G(x,y) f_n^+(y) dV(y)\right)dV(x).
\end{eqnarray}
Since $|G(x,y)-\dfrac{1}{2\pi}\log \dfrac{1}{|x-y|}|\leq C$ and $\left\|f_n\right\|_{L^1(M)}\leq C_1$, we get
\begin{eqnarray}
\label{regdecthm1etape1.3}
&&\int_{B_R (x_0)}\exp \left(\int_{B_{2R}(x_0)}G(x,y) f_n^+(y) dV(y)\right)dV(x)\nonumber\\
&\leq &C \sup_{y\in B_{2R}(x_0)}\left(\int_{B_R(x_0)}\left(\dfrac{1}{|x-y|}\right)^{\frac{\left\|f_n^{+}\right\|_{L^1(B_{2R} (x_0))}}{2\pi}}  \right)dV(x)\nonumber\\
&\leq & R^{2-\frac{\left\|f_n^{+}\right\|_{L^1(B_{2R} (x_0))}}{2\pi}}.
\end{eqnarray}
We note that, since $\left\|f_n^{+}\right\|_{L^1(B_{2R} (x_0))}\leq 2\pi$, we have $R^{2-\frac{\left\|f_n^{+}\right\|_{L^1(B_{2R} (x_0))}}{2\pi}}<+\infty$. Finally, thanks to (\ref{regdecthm1etape1.2}) and \eqref{regdecthm1etape1.3}, we get estimate \eqref{regdecthm1etape1.1}.
\begin{etape}
\label{regdecthm1etape2}
We have
\begin{eqnarray}
\label{regdecthm1entreetape}
&& \int_{B_R (x_0)} e^{v_n (x)}dV(x)\nonumber\\
 &\leq & C \exp \left(\displaystyle\sup_{x\in B_R(x_0)}\int_{M\backslash B_{2R}(x_0)} G(x,y) f_n^+ dV(y)\right)\nonumber \\
&\times & \exp \left(-\displaystyle\inf_{x\in B_r (y_0)}\int_{M \backslash B_r (y_0)}G(x,y) f_n^+ (y)dV(y)\right)\nonumber\\
&\times & R^{2-\frac{\left\|f_n^{+}\right\|_{L^1(B_{2R} (x_0))}}{2\pi}} r^{\frac{\left\| f_n^+\right\|_{L^1 (B_r (y_0))}}{2\pi}-2}\int_{B_r (y_0)} e^{v_n(x)}dV(x).
\end{eqnarray}
\end{etape}
\noindent\textit{Proof of Step 2.} Using the Green's representation formula, we have, for $x\in B_r(y_0)$,
\begin{eqnarray*}
v_n(x) &=& \bar{v}_n +\int_M G(x,y) f_n^+(y) dV(y) -\int_M G(x,y) \left(f_n^-(y) - Q_0(y) \right)dV(y).\nonumber \\
\end{eqnarray*}
Since, for all $q\geq 1$ and for all $x\in M$, $\left\|G(x,.)\right\|_{L^q(M)}\leq C_q$ where $C_q$ is a constant depending on $M$ and $q$, and since by hypothesis $\left\|f_n^-\right\|_{L^p(M)}\leq C_2$, using H\"{o}lder's inequality, we have
\begin{equation*}
\int_M G(x,y) f_n^- (y) dV(y) \leq \left(\int_M G(x,y)^{\frac{p}{p-1}}dV(y)\right)^{\frac{p-1}{p}} \left\|f_n^-\right\|_{L^p(M)}\leq C.
\end{equation*}
Thus, we get, for all $x\in B_r(y_0)$,
\begin{eqnarray*}
v_n(x)&\geq &\bar{v}_n +\int_M G(x,y) f_n^+(y) dV(y)-C \\
&\geq & -C +\bar{v}_n +\int_{B_r (y_0)}G(x,y) f_n^+ (y)dV(y)\\
&+& \inf_{x\in B_r (y_0)}\int_{M \backslash B_r (y_0)}G(x,y) f_n^+ (y)dV(y).
\end{eqnarray*}
Integrating $e^{v_n}$ on $B_r(y_0)$, one has
\begin{eqnarray}
\label{regdecthm1etape2.4}
\int_{B_r (y_0)} e^{v_n(x)}dV(x)&\geq & \exp \left(-C+\bar{v}_n+\displaystyle\inf_{x\in B_r (y_0)}\int_{M \backslash B_r (y_0)}G(x,y) f_n^+ (y)dV(y)\right)\nonumber\\
&\times & \int_{B_r (y_0)} \exp \left(\int_{B_r (y_0)}G(x,y) f_n^+ (y)dV(y)\right)dV(x)  .
\end{eqnarray}
Since, for $x\in B_r (y_0)$ and $y\in B_r (y_0)$, we have $\dfrac{1}{|x-y|}\geq \dfrac{1}{2r}$ and since $|G(x,y)-\dfrac{1}{2\pi}\log \dfrac{1}{|x-y|}|\leq C$, we get
\begin{eqnarray*}
&&\int_{B_r (y_0)} \exp \left(\int_{B_r (y_0)}G(x,y) f_n^+ (y)dV(y)\right)dV(x)\\
&\geq & C \int_{B_r (y_0)} \exp \left(\frac{1}{2\pi}\log \frac{1}{2r}\int_{B_r (y_0)}  f_n^+ (y)dV(y)\right)dV(x)\\
&\geq & C \left(\dfrac{1}{2r}\right)^{\frac{1}{2\pi}\int_{B_r (y_0)} f_n^+ (y) dV(y)} \int_{B_r (y_0)}dV(x)\\
&\geq & C r^{2-\frac{1}{2\pi}\int_{B_r (y_0)} f_n^+ (y) dV(y)}.
\end{eqnarray*}
Thanks to (\ref{regdecthm1etape2.4}) and the previous inequality, we derive that
\begin{eqnarray}
\label{regdecthm1etape2.1}
e^{\bar{v}_n} &\leq & C  \exp \left(-\displaystyle\inf_{x\in B_r (y_0)}\int_{M \backslash B_r (y_0)}G(x,y) f_n^+ (y)dV(y)\right)\nonumber\\
&\times &  r^{\frac{1}{2\pi}\int_{B_r (y_0)} f_n^+ (y) dV(y)-2}\int_{B_r (y_0)} e^{v_n(x)}dV(x).
\end{eqnarray}
Now, combining \eqref{regdecthm1etape1.1} and \eqref{regdecthm1etape2.1}, we obtain \eqref{regdecthm1entreetape}. 

\bigskip
\bigskip
\noindent\textit{Proof of Proposition \ref{regdecthm1}.} Thanks to \eqref{regdecthm1entreetape}, we see that, to prove the proposition, it is sufficient to prove that
\begin{eqnarray}
\label{11nov2012e2}
\sup_{x\in B_R(x_0)}\int_{M\backslash B_{2R}(x_0)}G(x,y)f_n^+ (y)dV(y) &-&\inf_{x\in B_r (y_0)}\int_{M\backslash B_{r}(y_0)}G(x,y)f_n^+ (y)dV(y)\nonumber\\
&\leq &C_\beta +\dfrac{1}{2\pi}\left(\log R\right) \int_{B_{2R}(x_0)}f_n^+ (y) dV(y)\nonumber\\
&-&\dfrac{1}{2\pi}\left(\log R\right) \int_{{B_r}(y_0)}f_n^+ (y) dV(y),
\end{eqnarray}
where $C_\beta$ is a positive constant depending on $\beta$. Let $x_1 \in B_R(x_0)$ and $x_2 \in B_r (y_0)$ be such that
\begin{eqnarray*}
&&\sup_{x\in B_R(x_0)}\int_{M\backslash B_{2R}(x_0)}G(x,y)f_n^+ (y)dV(y) -\inf_{x\in B_r (y_0)}\int_{M\backslash B_{r}(y_0)}G(x,y)f_n^+ (y)dV(y)\\
&=& \int_{M\backslash B_{2R}(x_0)}G(x_1,y)f_n^+ (y)dV(y) -\int_{M\backslash B_{r}(y_0)}G(x_2,y)f_n^+ (y)dV(y).
\end{eqnarray*}
Since $B_r (y_0) \subset B_{4(\beta +1)R}(x_0)$ and $B_{2R}(x_0)\subset B_{4(\beta +1)R}(x_0)$, we can write
$$\int_{M\backslash B_{2R}(x_0)}G(x_1,y)f_n^+ (y)dV(y) -\int_{M\backslash B_{r}(y_0)}G(x_2,y)f_n^+ (y)dV(y)=I+II+III,$$
where
$$I=\int_{M\backslash B_{4(\beta +1)R}(x_0)}G(x_1,y) f_n^+(y)dV(y)-\int_{M\backslash B_{4(\beta +1)R}(x_0)}G(x_2,y) f_n^+(y)dV(y), $$
$$II=\int_{ B_{4(\beta +1)R}(x_0)\backslash B_{2R}(x_0)}G(x_1,y) f_n^+(y)dV(y)$$
and
$$III=-\int_{ B_{4(\beta +1)R}(x_0)\backslash B_{r}(y_0)}G(x_2,y) f_n^+(y)dV(y).$$
We first estimate $I$. For $y\in M\backslash B_{4(\beta +1)R}(x_0)$, one can easily check that there exists a constant $C_\beta$ depending on $\beta$ such that
\begin{equation*}
\label{21juin2012e10}
C_\beta \leq \dfrac{|y-x_2|}{|y-x_1|}\leq \dfrac{1}{C_\beta}.
\end{equation*}
This implies that there exists a constant $\tilde{C}_\beta$ such that, for all $y\in M\backslash B_{4(\beta +1)R}(x_0)$,
$$|G(x_1,y)-G(x_2,y)|\leq C+\dfrac{1}{2\pi}\left|\log \dfrac{|x_2-y|}{|x_1-y|}\right|\leq \tilde{C}_\beta.$$
Therefore, since $\left\|f_n\right\|_{L^1(M)}\leq C_1$, we have
\begin{eqnarray}
\label{regdecthm1etape3.3}
I&=& \int_{M\backslash B_{4(\beta +1)R}(x_0)}G(x_1,y) f_n^+(y)dV(y)- \int_{M\backslash B_{4(\beta +1)R}(x_0)}G(x_2,y) f_n^+(y)dV(y)\nonumber\\
& \leq & \tilde{C}_\beta \left\|f_n\right\|_{L^1(M)}\leq \tilde{C}_\beta C_1.
\end{eqnarray} 
Let us now estimate $II$. Since $x_1 \in B_R (x_0)$, for $y\in B_{4(\beta +1)R}(x_0)\backslash B_{2R}(x_0)$, we note that $|x_1-y|\geq R$. Thus recalling the asymptotic of the Green's function and that $\left\|f_n\right\|_{L^1(M)}\leq C_1$, we deduce that 
\begin{eqnarray}
\label{regdecthm1etape3.4}
II&=&\int_{ B_{4(\beta +1)R}(x_0)\backslash B_{2R}(x_0)}G(x_1,y) f_n^+(y)dV(y)\nonumber\\
 & \leq & \dfrac{1}{2\pi}\left(\log \dfrac{1}{R}\right)\int_{ B_{4(\beta +1)R}(x_0)\backslash B_{2R}(x_0)} f_n^+(y) dV(y)+C C_1.
\end{eqnarray} 
Finally, noting that for $x_2 \in B_r(y_0)$ and $y\in  B_{4(\beta +1)R}(x_0)\backslash B_{r}(y_0)$, we have $|x_2-y|\leq 5(\beta +1)R$, we obtain
\begin{eqnarray}
\label{regdecthm1etape3.5}
III&=&-\int_{ B_{4(\beta +1)R}(x_0)\backslash B_{r}(y_0)}G(x_2,y) f_n^+(y)dV(y)\nonumber\\
 & \leq & \dfrac{1}{2\pi}\left(\log R\right)\int_{ B_{4(\beta +1)R}(x_0)\backslash B_{r}(y_0)} f_n^+(y) dV(y)+C.
\end{eqnarray} 
Therefore, from (\ref{regdecthm1etape3.3}), (\ref{regdecthm1etape3.4}) and (\ref{regdecthm1etape3.5}), we get
\begin{eqnarray*}
&&\int_{M\backslash B_{2R}(x_0)}G(x_1,y)f_n^+ (y)dV(y) -\int_{M\backslash B_{r}(y_0)}G(x_2,y)f_n^+ (y)dV(y)\\
&\leq &C_\beta -\dfrac{1}{2\pi}\left(\log R\right)\int_{ B_{4(\beta +1)R}(x_0)\backslash B_{2R}(x_0)} f_n^+(y) dV(y)\\
&+& \dfrac{1}{2\pi}\left(\log R\right)\int_{ B_{4(\beta +1)R}(x_0)\backslash B_{r}(y_0)} f_n^+(y) dV(y) \\
&\leq &C_\beta +\dfrac{1}{2\pi}(\log R )\int_{B_{2R}(x_0)}f_n^+ (y) dV(y)\\
&-&\dfrac{1}{2\pi}(\log R )\int_{{B_r}(y_0)}f_n^+ (y) dV(y).
\end{eqnarray*}
Hence \eqref{11nov2012e2} holds. This achieves the proof of Proposition \ref{regdecthm1}.

\end{proof}

\subsection{Proof of Theorem \ref{regdecthm3}.}
To prove Theorem \ref{regdecthm3}, we will need three lemma. We first prove that the integral of $e^{v_n}$ on some annulus, centered in a blow-up point, tends to $0$ at infinity. 
\begin{lem}
\label{regdecprop1}
Let $(x_n,r_n)_n$ be a blow-up for $(v_n)_n$, and $(R_n)_n$, $(S_n)_n$ two sequences of positive numbers satisfying the following properties :
\begin{equation}
\label{regdecprop1.1}
0<2 R_n \leq \dfrac{S_n}{4},
\end{equation}
\begin{equation}
\label{regdecprop1.2}
 \dfrac{R_n}{r_n}\underset{n\rightarrow +\infty}{\longrightarrow}+\infty,
\end{equation}
and there exists $\bar{C}>0$ (not depending on $n$) such that
\begin{equation}
\label{regdecprop1.3}
\forall B_s(y)\subset B_{S_n}(x_n)\backslash B_{R_n}(x_n),\ if \int_{B_s (y)}e^{v_n}dV \geq \dfrac{\pi}{\rho}\ then\ s\geq \bar{C}|y-x_n|.
\end{equation}
Then, we have
\begin{equation}
\label{regdecprop1.4}
\lim_{n\rightarrow +\infty}\int_{B_{\frac{S_n}{2}}(x_n)\backslash B_{2R_n}(x_n)}e^{v_n}dV= 0.
\end{equation}
\end{lem}

\begin{proof}
Since $(x_n,r_n)$ is a blow-up for $(v_n)_n$, then, from Proposition \ref{mp2}, we can choose $(b_n)_n$ such that
\begin{equation}
\label{regdecprop1.6}
\dfrac{b_n r_n}{R_n}\underset{n\rightarrow +\infty}{\longrightarrow} 0,
\end{equation}
and
\begin{equation}
\label{regdecprop1.5}
\lim_{n\rightarrow +\infty}\int_{B_{b_n r_n}(x_n)} e^{v_n}dV= \dfrac{8\pi}{\rho}.
\end{equation}
We will divide the proof into two steps. In the first one, we will obtain an upper bound of $\displaystyle\int_{B_{2R}(x_n)\backslash B_{R}(x_n)}e^{v_n}dV$ for $2R_n \leq R \leq \dfrac{S_n}{4}$. In the second step, we will establish (\ref{regdecprop1.4}) by dividing the annulus $B_{\frac{S_n}{2}}(x_n)\backslash B_{2R_n}(x_n)$ into smaller annuli on which the estimate of the first step will hold.
\setcounter{etape}{0}
\begin{etape}
\label{regdecprop1etape1}
Let $R>0$ be such that $2R_n \leq R \leq \dfrac{S_n}{4}$, then we have
$$\int_{B_{2R}(x_n)\backslash B_{R}(x_n)}e^{v_n}dV \leq C \left(\dfrac{b_n r_n}{R}\right)^{2+o_n(1)},$$
where $o_n(1)\underset{n\rightarrow +\infty}{\longrightarrow}0 $.
\end{etape}
\noindent\textit{Proof of Step 1.} First, for $y\in B_{2R}(x_n)\backslash B_{R}(x_n)$, we have
$$B_{\frac{\bar{C}R}{8}}(y) \subset B_{S_n}(x_n)\backslash B_{R_n}(x_n),$$
where $\bar{C}$ is the constant defined in (\ref{regdecprop1.3}) (without loss of generality, we can assume that $\bar{C}<1$). Thus we deduce from (\ref{regdecprop1.3}) that, for $y\in B_{2R}(x_n)\backslash B_{R}(x_n)$,
\begin{equation}
\label{regdecprop1etape1.1}
\int_{B_{\frac{\bar{C}R}{8}}(y)} e^{v_n}dV < \dfrac{\pi}{\rho}.
\end{equation}
On the other hand, $B_{2R}(x_n)\backslash B_{R}(x_n)$ can be covered by a finite number (not depending on $n$) of balls $B_{\frac{\bar{C}R}{16}}(y_i)$ where $y_i\in B_{2R}(x_n)\backslash B_{R}(x_n)$. We will use Proposition \ref{regdecthm1} for $f_n=\rho e^{v_n}+h_n e^{v_n}$, $x_0=y_i$, $y_0=x_n$, $r=b_n r_n$ and $R$ (the one defined in Proposition \ref{regdecthm1}) substituted by $\tilde{R}=\dfrac{\bar{C}R}{16}$. Let's show that the hypothesis of the proposition hold true. Since $\displaystyle\int_M e^{v_n}dV =1$, $\forall n\geq 0$, we have, using Holder's inequality and (\ref{regdeceq1})(i), that
$$\left\|f_n \right\|_{L^1(M)}\leq \rho \int_M e^{v_n}dV +\int_M |h_n|e^{v_n}dV\leq \rho + \left(\int_M |h_n|^2 e^{v_n}dV\right)^{\frac{1}{2}}\leq C.$$
Since $y_i \in A_{R,2R}(x_n)$, there exists a constant $\beta >0$ such that $$|y_i-x_n|\leq 2R\leq \beta\tilde{R}.$$ 
It is easy to check, using \eqref{regdecprop1.6}, that 
$$r=b_n r_n\leq \tilde{R}=\dfrac{\bar{C}R}{16}.$$
Therefore, it only remains to check that there exists $p>1$ such that 
$$\left\|f_n^- \right\|_{L^p(M)}\leq C,$$
and
$$\int_{B_{\frac{\bar{C}R}{8}}(y_i)}f_n^+dV \leq 2\pi.$$
The first inequality is a direct consequence of (\ref{regdeceq1})(ii). The second one follows from (\ref{regdecprop1etape1.1}) and (\ref{regdeceq1})(i) which imply, using Holder's inequality, that
\begin{eqnarray*}
\int_{B_{\frac{\bar{C}R}{8}}(y_i) }f_n^+ dV &\leq & \rho \int_{B_{\frac{\bar{C}R}{8}}(y_i)} e^{v_n}dV+ \left(\int_M h_n^2 e^{v_n}dV\right)^{\frac{1}{2}}\\
& \leq & \pi+o_n(1)\leq 2\pi .
\end{eqnarray*}
Thus, Proposition \ref{regdecthm1} gives us that
\begin{eqnarray*}
\int_{B_{\frac{\bar{C}R}{16}}(y_i)}e^{v_n}dV&\leq &C \left(\dfrac{b_n r_n}{R}\right)^{\frac{\left\|f_n^+\right\|_{L^1 (B_{r_n b_n}(x_n))}}{2\pi}-2}\int_{B_{b_n r_n}(x_n) } e^{v_n}dV\\
&\leq &C \left(\dfrac{b_n r_n}{R}\right)^{\frac{\left\|f_n^+\right\|_{L^1 (B_{r_n b_n}(x_n))}}{2\pi}-2}.
\end{eqnarray*}
From (\ref{regdecprop1.5}), since $\displaystyle\lim_{n\rightarrow +\infty}\int_{B_{b_n r_n}(x_n)} |h_n| e^{v_n}dV=0$ and $\rho e^{v_n} \leq f_n^+ \leq \rho e^{v_n}+|h_n|e^{v_n}$, we get 
$$\frac{\left\|f_n^+\right\|_{L^1(B_{r_n b_n}(x_n))}}{2\pi}-2=2+o_n(1).$$
Therefore, we obtain
$$\int_{B_{\frac{\bar{C}R}{16}}(y_i)}e^{v_n}dV\leq C \left(\dfrac{b_n r_n}{R}\right)^{2+o_n(1)}.$$
Since $B_{2R}(x_n)\backslash B_{R}(x_n)$ is covered by a finite number not depending on $n$ of balls $B_{\frac{\bar{C}R}{16}}(y_i)$, we finally derive that
\begin{equation}
\label{12nov2012e2}
\int_{B_{2R}(x_n)\backslash B_{R}(x_n)}e^{v_n}dV \leq C \left(\dfrac{b_n r_n}{R}\right)^{2+o_n(1)}.
\end{equation}
\begin{etape}
We have
$$\underset{n\rightarrow +\infty}{\lim}\int_{B_{\frac{S_n}{2}}(x_n)\backslash B_{2R_n}(x_n)} e^{v_n}dV = 0.$$
\end{etape}
\noindent\textit{Proof of Step 2.} We apply Step \ref{regdecprop1etape1} to $\tilde{R}_j=2^j (2 R_n)$ , $0\leq j\leq j_n$ with $j_n$ such that $2^{j_n}\in \left[\dfrac{S_n}{16 R_n},\dfrac{S_n}{8 R_n}  \right]$. By definition, we see that $\displaystyle \bigcup_{j=0}^{j_n} B_{2\tilde{R}_j}(x_n)\backslash B_{\tilde{R}_j}(x_n)$ covers $B_{\frac{S_n}{4}}(x_n)\backslash B_{2R_n}(x_n)$. We deduce that
$$\int_{B_{\frac{S_n}{2}}(x_n)\backslash B_{2R_n}(x_n)} e^{v_n}dV \leq \sum_{j=0}^{j_n}\int_{B_{2\tilde{R}_j}(x_n)\backslash B_{\tilde{R}_j}(x_n) }e^{v_n}dV + \int_{B_{\frac{S_n}{2}}(x_n)\backslash B_{\frac{S_n}{4}}(x_n)}e^{v_n}dV.$$ 
From Step \ref{regdecprop1etape1}, since $2R_n \leq \tilde{R}_j\leq \dfrac{S_n}{4}$, for all $0\leq j\leq j_n$,
$$\int_{B_{\frac{S_n}{2}}(x_n)\backslash B_{2R_n}(x_n)} e^{v_n}dV \leq C\sum_{j=0}^{j_n}\left(\dfrac{b_n r_n}{2^j 2R_n}\right)^{2+o_n(1)} +C\left( \dfrac{b_nr_n} {\frac{S_n}{4}}\right)^{2+o_n(1)}, $$
and since by definition , $\dfrac{S_n}{4}\geq 2R_n$, we deduce that
\begin{eqnarray*}
\int_{B_{\frac{S_n}{2}}(x_n)\backslash B_{2R_n}(x_n)} e^{v_n}dV &\leq &C \left(\dfrac{b_n r_n}{R_n}\right)^{2+o_n(1)}\left(\sum_{j=0}^{\infty}2^{-2j} +1\right)\leq  C \left(\dfrac{b_n r_n}{R_n}\right)^{2+o_n(1)}.
\end{eqnarray*}
From (\ref{regdecprop1.6}), $\dfrac{b_n r_n}{R_n}\underset{n\rightarrow +\infty}{\longrightarrow} 0$, we finally arrive at
$$\underset{n\rightarrow +\infty}{\lim}\int_{B_{\frac{S_n}{2}}(x_n)\backslash B_{2R_n}(x_n)} e^{v_n}dV = 0.$$
Therefore, Lemma \ref{regdecprop1} is established.
\end{proof}

Now, we will prove that, assuming there exists $k$ blow-up for $(v_n)_n$, there is no volume on annuli centered in these points. 
\begin{lem}
\label{regdecprop2}
Let $(x_n^1,r_n^1)_n,\ldots,(x_n^k,r_n^k)_n$ be $k$ blow-up for $(v_n)_n$ with $k\geq 2$, and suppose there exist sequences $(R_n^1)_n,\ldots,(R_n^k)_n$ satisfying the following properties :\newline
1)\begin{equation}
\label{regdecprop2.1.1} 
\dfrac{R_n^i}{r_n^i}\underset{n\rightarrow +\infty}{\longrightarrow}+\infty, \ \forall i\in \{1,\ldots ,k\}. 
\end{equation}
2)
\begin{equation}
\label{regdecprop2.1.2}
 \dfrac{R_n^i}{\displaystyle\inf_{i\neq j\in \{1,\ldots,k\}}|x_n^i -x_n^j|}\underset{n\rightarrow +\infty}{\longrightarrow} 0,\ \forall i\in \{1,\ldots ,k\}.
\end{equation} 
3) If $D_n =\displaystyle\max_{i\neq j\in \{1,\ldots,k\}}|x_n^i -x_n^j|$, we suppose that there exists a constant $\bar{C}>0$ such that
\begin{eqnarray}
\label{regdecprop2.1.3}
&&\forall B_s (y)\subset \bigcup_{i=1}^k B_{4D_n}(x_n^i)\backslash \bigcup_{i=1}^k B_{R_n^i}(x_n^i),\ if\ \int_{B_s (y)}e^{v_n}dV \geq \dfrac{\pi}{\rho} \nonumber \\
&&then \ s\geq \bar{C} d_n (y)=\bar{C} \displaystyle\inf_{i\in \{1,\ldots , k\}}|y-x_n^i|.
\end{eqnarray}
Then we have
$$\lim_{n\rightarrow +\infty} \int_{\bigcup_{i=1}^k B_{2D_n}(x_n^i)\backslash \bigcup_{i=1}^k B_{2R_n^i}(x_n^i)}e^{v_n}dV=0.$$
\end{lem}

\begin{proof}
We will proceed by induction on $k$. For $k=2$, we have $D_n=|x_n^1 -x_n^2|$. From (\ref{regdecprop2.1.1}) and (\ref{regdecprop2.1.2}), we get
$$\dfrac{R_n^1}{r_n^1}\underset{n\rightarrow +\infty}{\longrightarrow} +\infty,$$
and
$$\dfrac{R_n^1}{D_n}\underset{n\rightarrow +\infty}{\longrightarrow} 0.$$
Moreover, if $B_s (y)\subset B_{\frac{D_n}{2}}(x_n^1)\backslash B_{R_n^1}(x_n^1)\subset \displaystyle\bigcup_{i=1}^2 B_{4D_n}(x_n^i)\backslash \bigcup_{i=1}^2 B_{R_n^i}(x_n^i)$ then $d_n (y) =|y-x_n^1|$. Therefore, using (\ref{regdecprop2.1.3}), Lemma \ref{regdecprop1} for $x_n=x_n^1$ with $R_n =R_n^1$ and $S_n =\dfrac{D_n}{2}$, implies that
\begin{equation}
\label{regdecprop2.1}
\lim_{n\rightarrow +\infty}\int_{ B_{\frac{D_n}{4}}(x_n^1)\backslash  B_{2R_n^1}(x_n^1)}e^{v_n}dV= 0.
\end{equation}
In the same way, we also have
\begin{equation}
\label{regdecprop2.2}
\lim_{n\rightarrow +\infty}\int_{ B_{\frac{D_n}{4}}(x_n^2)\backslash  B_{2R_n^2}(x_n^2)}e^{v_n}dV=0.
\end{equation}
We can cover  $\displaystyle\bigcup_{i=1}^2 B_{2D_n}(x_n^i) \backslash \bigcup_{i=1}^2 B_{\frac{D_n}{4}}(x_n^i)$ with a finite number $N$ (not depending on $n$) of balls $B_{\frac{\bar{C}D_n}{16}}(y_j)$, where $\bar{C}$ is the constant defined in \eqref{regdecprop2.1.3} and $y_j\in \displaystyle\bigcup_{i=1}^2 B_{2D_n}(x_n^i) \backslash \bigcup_{i=1}^2 B_{\frac{D_n}{4}}(x_n^i)$, such that
$$B_{\frac{\bar{C}D_n}{16}}(y_j) \subset \displaystyle\bigcup_{i=1}^2 B_{4D_n}(x_n^i) \backslash \displaystyle\bigcup_{i=1}^2 B_{R_n^i}(x_n^i).$$
It is easy to check, by using \eqref{regdecprop2.1.3}, that 
$$\int_{B_{\frac{\bar{C}D_n}{8}}(y_j)}e^{v_n}dV \leq \dfrac{\pi}{\rho}.$$
Now following the same argument as to get \eqref{12nov2012e2}, we obtain
\begin{equation}
\label{24maie4}
\int_{\bigcup_{i=1}^2 B_{2D_n}(x_n^i) \backslash \bigcup_{i=1}^2 B_{\frac{D_n}{4}}(x_n^i)}e^{v_n}dV \leq C \left(\dfrac{b_n r_n^1}{D_n}\right)^{2+o_n(1)},
\end{equation}
 where $b_n \underset{n\rightarrow +\infty}{\longrightarrow} +\infty$ satisfies $\dfrac{b_n r_n^1}{R_n^1}\underset{n\rightarrow +\infty}{\longrightarrow} 0$. Thus, since $\dfrac{b_n r_n^1}{D_n}= \dfrac{b_n r_n^1}{R_n^1}\dfrac{R_n^1}{D_n}\underset{n\rightarrow +\infty}{\longrightarrow} 0$, we deduce that
\begin{equation}
\label{regdecprop2.3}
\lim_{n\rightarrow +\infty}\int_{\bigcup_{i=1}^2 B_{2D_n}(x_n^i) \backslash \bigcup_{i=1}^2 B_{\frac{D_n}{4}}(x_n^i)} e^{v_n}dV= 0.
\end{equation}
So, from (\ref{regdecprop2.1}), (\ref{regdecprop2.2}) and (\ref{regdecprop2.3}), we find
$$\lim_{n\rightarrow +\infty}\int_{\bigcup_{i=1}^2 B_{2D_n}(x_n^i)\backslash \bigcup_{i=1}^2 B_{2R_n^i}(x_n^i)}e^{v_n}dV= 0.$$
Therefore the lemma holds for $k=2$. Now suppose that the lemma holds for $k$. Let $d_n=\displaystyle\inf_{i\neq j \in \{1,\ldots ,k+1\}}|x_n^i-x_n^{j}|$. We will consider two cases depending on the value of $\dfrac{d_n }{D_n}$ when $n\rightarrow +\infty$.
\newline
\noindent \textit{First case : $\displaystyle \dfrac{d_n }{D_n}\underset{n\rightarrow +\infty}{\not\longrightarrow} 0$.}\newline
In this case, there exists a constant $C>0$ such that $d_n \geq CD_n$. Thus, we have 
\begin{equation}
\label{5avril20122secondcase1}
C\leq \dfrac{d_n }{D_n}\leq 1.
\end{equation}
Using \eqref{5avril20122secondcase1}, we can cover  $\displaystyle\bigcup_{i=1}^{k+1} B_{2D_n}(x_n^i) \backslash \bigcup_{i=1}^{k+1} B_{\frac{D_n}{4}}(x_n^i)$ with a finite number $N$ (not depending on $n$) of balls $B_{\frac{\bar{C}D_n}{16}}(y_j)$, where $\bar{C}$ is the constant defined in \eqref{regdecprop2.1.3} and $y_j\in \displaystyle\bigcup_{i=1}^{k+1} B_{2D_n}(x_n^i) \backslash \bigcup_{i=1}^{k+1} B_{\frac{D_n}{4}}(x_n^i)$, such that
$$B_{\frac{\bar{C}D_n}{16}}(y_j) \subset \displaystyle\bigcup_{i=1}^{k+1} B_{4D_n}(x_n^i) \backslash \displaystyle\bigcup_{i=1}^{k+1} B_{R_n^i}(x_n^i).$$
Following the same arguments as the one used for the case $k=2$, we obtain that
$$\lim_{n\rightarrow +\infty}\int_{\bigcup_{i=1}^{k+1} B_{2D_n}(x_n^i)\backslash \bigcup_{i=1}^{k+1} B_{2R_n^i}(x_n^i)}e^{v_n}dV= 0.$$
 Let's now consider the second case.\\
\textit{Second case : $\dfrac{d_n }{D_n}\underset{n\rightarrow +\infty}{\longrightarrow} 0.$}\newline
We can assume, up to relabelling sequences $(x_n^i)_{n}$, that $d_n=|x_n^1 -x_n^{k+1}|$. We define the set $X_1$ of blow-up points for which the distance to $x_n^1$ is comparable to $d_n$, that is to say
$$X_1 = \left\{(x_n^j)_{n } :\exists C_j >0 \ such\ that\ |x_n^j -x_n^1|\leq C_j d_n,\ \forall n\geq 0\right\}.$$
Notice that $(x_n^1)_{n}$ and $(x_n^{k+1})_{n}$ are in $X_1$ and that $card (X_1)< k+1$ since we assume $\dfrac{d_n }{D_n}\underset{n\rightarrow +\infty}{\longrightarrow} 0$. Up to relabelling sequences $(x_n^j)_{n}$, for $j\neq 1,k+1$, we can assume that
$$X_1= \left\{(x_n^j)_{n}, (x_n^{k+1})_{n} \ :\ j\in \{1,\ldots ,l_0+1\}  \right\},$$
with $0\leq l_0\leq k-2$. Let $C= \max \left\{C_j\ :\ j\in \{2,\ldots , l_0+1 ,k+1 \right\}$. We have
\begin{equation}
\label{regdecprop2.1bis}
|x_n^1 -x_n^j|\leq Cd_n, \ \forall j\in \{2,\ldots , l_0+1 ,k+1\},
\end{equation}
and
\begin{equation}
\label{regdecprop2.2bis}
\dfrac{|x_n^1-x_n^j|}{d_n}\underset{n\rightarrow +\infty}{\longrightarrow} +\infty, \ \forall j\in \{l_0+2,\ldots ,k\}.
\end{equation}
One can check that the induction hypothesis of the lemma holds for the sequences
$$(x_n^1,r_n^1)_{n}, (x_n^{l_0+2},r_n^{l_0+2})_{n},\ldots , (x_n^k,r_n^k)_{n}$$
associated to resp. $(2Cd_n)_n, (R_n^{l_0+2})_n,\ldots ,(R_n^k)_n$. It follows that
\begin{equation}
\label{regdecprop2.4}
I=\displaystyle\lim_{n\rightarrow +\infty}\int_{\bigcup_{i\in \{1,l_0+2,\ldots ,k\}} B_{2\tilde{D}_n}(x_n^i)\backslash \left(\bigcup_{i=l_0+2}^k B_{2R_n^i}(x_n^i)\cup B_{4Cd_n}(x_n^1)\right)}e^{v_n}dV= 0,
\end{equation}
where $\tilde{D}_n=\displaystyle\max_{i\neq j \in \{1,l_0+2,\ldots , k\}}|x_n^i -x_n^j|$. Let us show that
$$\displaystyle\lim_{n\rightarrow +\infty}\int_{ B_{4Cd_n}(x_n^1)\backslash \bigcup_{i\in \{1,\ldots ,l_0+1,k+1\}} B_{2R_n^i}(x_n^i)}e^{v_n}dV= 0.$$
One can check that the hypothesis to apply Lemma \ref{regdecprop1} for  $x_n^i$, $i\in \{1, \ldots , l_0+1,k+1\}$, with $R_n =R_n^i$, $S_n=\dfrac{d_n}{2}$ hold. Thus, we have, for all $i\in \{1, \ldots , l_0+1, k+1\}$,
\begin{equation}
\label{regdecprop2.5}
II_i=\displaystyle\lim_{n\rightarrow +\infty}\int_{A_{2R_n^i ,\frac{d_n}{4}}(x_n^i)}e^{v_n}dV= 0.
\end{equation}
On the other hand, proceeding as in the case $k=2$, we can cover the set
$$B_{4Cd_n}(x_n^1) \backslash \bigcup_{i\in \{1,\ldots , l_0+1, k+1\}} B_{\frac{d_n}{4}}(x_n^i)$$ by a finite number $N$ not depending on $n$ of balls  
such that
$$B_{\frac{\bar{C}d_n}{16}}(y_j) \subset \bigcup_{i=1}^{k+1} B_{4D_n}(x_n^i) \backslash \bigcup_{i=1}^{k+1} B_{R_n^i}(x_n^i),$$
and
$$\int_{B_{\frac{\bar{C} d_n}{8}}(y_j)}e^{v_n}dV \leq \dfrac{\pi}{\rho}.$$
Therefore, we obtain, in a similar way as for \eqref{24maie4}, that there exists a sequence $(b_n)_n$, $b_n \underset{n\rightarrow +\infty}{\longrightarrow}+\infty$ arbitrarily slowly, such that 
\begin{equation}
\label{25maie5}
\dfrac{b_n r_n^1}{R_n^1}\underset{n\rightarrow +\infty}{\longrightarrow}0,
\end{equation}
and
$$\int_{B_{\frac{\bar{C}d_n}{8}}(y_j) }e^{v_n}dV \leq C\left(\dfrac{b_n r_n^1}{d_n} \right)^{2+o_n(1)}.$$
From \eqref{regdecprop2.1.2} and \eqref{25maie5}, we have $\dfrac{b_n r_n^1}{d_n}\underset{n\rightarrow +\infty}{\longrightarrow}0$. Hence, since $$\displaystyle B_{4Cd_n}(x_n^1) \backslash \bigcup_{i\in \{1,\ldots , l_0+1, k+1\}} B_{\frac{d_n}{4}}(x_n^i)$$ is covered by a finite number $N$ not depending on $n$ of balls $B_{\frac{\bar{C}d_n}{16}}(y_j)$, we obtain that
\begin{equation}
\label{regdecprop2.6}
III=\displaystyle\lim_{n\rightarrow +\infty}\int_{ B_{4Cd_n}(x_n^1) \backslash \bigcup_{i\in \{1,\ldots , l_0+1, k+1\}} B_{\frac{d_n}{4}}(x_n^i)  } e^{v_n}dV= 0.
\end{equation}
Therefore, (\ref{regdecprop2.4}), (\ref{regdecprop2.5}) and (\ref{regdecprop2.6}) yield to
\begin{eqnarray}
\label{avril2012avantderlem2}
&&\displaystyle\lim_{n\rightarrow +\infty}\int_{\bigcup_{i\in \{1,l_0 +2,\ldots , k\}} B_{2\tilde{D}_n}(x_n^i) \backslash \bigcup_{i=1}^{k+1} B_{2R_n^i}(x_n^i)}e^{v_n} dV\nonumber \\
& =& I+\displaystyle\sum_{i\in \{1,\ldots ,l_0+1,k+1\}}II_i+III=  0.
\end{eqnarray}
Now, we claim that 
\begin{equation}
\label{25maie8}
\displaystyle\lim_{n\rightarrow +\infty} \int_{\bigcup_{i=1}^{k+1} B_{2D_n}(x_n^i)\backslash \bigcup_{i=1}^{k+1} B_{2R_n^i}(x_n^i)}e^{v_n}dV=0.
\end{equation}
It is not difficult to check that
$$\displaystyle\bigcup_{i=1}^{k+1} B_{\frac{3}{2}D_n}(x_n^i) \subset \displaystyle\bigcup_{i\in \{1,l_0 +2,\ldots ,k\}} B_{2\tilde{D}_n}(x_n^i).$$
Therefore we deduce, from \eqref{avril2012avantderlem2}, that
\begin{equation}
\label{6avril2012eq1}
\displaystyle\lim_{n\rightarrow +\infty}\int_{\bigcup_{i=1}^{k+1}B_{\frac{3}{2}D_n}(x_n^i) \backslash \bigcup_{i=1}^{k+1} B_{2R_n^i}(x_n^i)  } e^{v_n}dV = 0.
\end{equation}
To prove \eqref{25maie8}, it is sufficient, from \eqref{6avril2012eq1}, to show that
$$\displaystyle\lim_{n\rightarrow +\infty}\int_{\bigcup_{i=1}^{k+1}B_{2D_n}(x_n^i) \backslash \bigcup_{i=1}^{k+1}B_{\frac{3}{2}D_n}(x_n^i)}e^{v_n}dV =0.$$
In the same way as for the case $k=2$, we cover the set
$$\bigcup_{i=1}^{k+1}B_{2D_n}(x_n^i) \backslash \bigcup_{i=1}^{k+1}B_{\frac{3}{2}D_n}(x_n^i)$$ by a finite number $N$ (not depending on $n$) of balls $B_{\frac{\bar{C}D_n}{16}}(y_j)$ ,
$$y_j \in \bigcup_{i=1}^{k+1}B_{2D_n}(x_n^i) \backslash \bigcup_{i=1}^{k+1}B_{\frac{3}{2}D_n}(x_n^i),$$
such that
$$\displaystyle\lim_{n\rightarrow +\infty}\int_{B_{\frac{\bar{C}D_n}{8}}(y_j) }e^{v_n}dV = 0.$$
Thus, we obtain 
\begin{equation}
\label{avril2012avantderlem1}
\displaystyle\lim_{n\rightarrow +\infty}\int_{\bigcup_{i=1}^{k+1}B_{2D_n}(x_n^i) \backslash \bigcup_{i=1}^{k+1}B_{\frac{3}{2}D_n}(x_n^i)}e^{v_n}dV = 0.
\end{equation}
Finally, from \eqref{6avril2012eq1} and \eqref{avril2012avantderlem1}, we get
$$\displaystyle\lim_{n\rightarrow +\infty}\int_{\bigcup_{i=1}^{k+1}B_{2D_n}(x_n^i) \backslash \bigcup_{i=1}^{k+1} B_{2R_n^i}(x_n^i)  } e^{v_n}dV= 0.$$
This achieves the proof of the lemma.
\end{proof}
Finally, we have :
\begin{lem}
\label{regdecprop4maii}
\label{regdecprop2ter}
Let $(x_n^1,r_n^1)_n,\ldots,(x_n^k,r_n^k)_n$ be $k$ blow-up for $(v_n)_n$ with $k\geq 2$, and suppose that there exist $k$ sequences $(R_n^1)_n,\ldots,(R_n^k)_n$ satisfying the following properties :\newline
1)\begin{equation}
\label{regdecprop2.1.1ter} 
\dfrac{R_n^i}{r_n^i}\underset{n\rightarrow +\infty}{\longrightarrow}+\infty, \ \forall i\in \{1,\ldots ,k\}. 
\end{equation}
2)
\begin{equation}
\label{regdecprop2.1.2ter}
 \dfrac{R_n^i}{\displaystyle\inf_{i\neq j\in \{1,\ldots,k\}}|x_n^i -x_n^j|}\underset{n\rightarrow +\infty}{\longrightarrow} 0,\ \forall i\in \{1,\ldots ,k\}.
\end{equation} 
3) There exists $\bar{C}>0$ such that
\begin{eqnarray}
\label{regdecprop2.1.3ter}
&&\forall B_s (y)\subset M\backslash \bigcup_{i=1}^k B_{R_n^i}(x_n^i),\ if\ \int_{B_s (y)}e^{v_n}dV \geq \dfrac{\pi}{\rho} \nonumber \\
&&then \ s\geq \bar{C} d_n (y),\ where\ d_n(y)=\displaystyle\inf_{i\in \{1,\ldots , k\}}|y-x_n^i|.
\end{eqnarray}
Then
$$\displaystyle\lim_{n\rightarrow +\infty} \int_{ M\backslash \bigcup_{i=1}^k B_{2R_n^i}(x_n^i)}e^{v_n}dV=0.$$
\end{lem}
\begin{proof}
The proof relies on the same technics as the ones used in the previous lemma, therefore we omit it.
\end{proof}
\bigskip
\vspace{12pt}

\noindent\textit{Proof of Theorem \ref{regdecthm3}.} Suppose that $(v_n)_n\subseteq H^2(M)$ is not compact, we claim that there exist $k$ blow-up $$(x_n^1, r_n^1)_{n},\ldots ,(x_n^k, r_n^k)_{n}$$ with $1\leq k\leq [\dfrac{\rho}{8\pi}]$ (where $[\dfrac{\rho}{8\pi}]$ stands for the integer part of $\dfrac{\rho}{8\pi}$), and there exist $k$ sequences $(R_n^1)_n,\ldots ,(R_n^k)_n >0$ satisfying the following properties :\newline
1)
\begin{equation}
\label{regdecprop2.1.128mai} 
\dfrac{R_n^i}{r_n^i}\underset{n\rightarrow +\infty}{\longrightarrow}+\infty, \ \forall i\in \{1,\ldots ,k\}. 
\end{equation}
2)
\begin{equation}
\label{regdecprop3.1}
 \lim_{n\rightarrow +\infty} \int_{B_{2R_n^i}(x_n^i)} e^{v_n}dV=\dfrac{8\pi}{\rho},\ \forall i\in \{1,\ldots , k\}.
 \end{equation}
3)
\begin{eqnarray}
\label{regdecprop3.2} 
&& \dfrac{R_n^i}{\underset{ i\neq j\in \{1,\ldots ,k\} }{\inf} |x_n^i-x_n^j|}\underset{n\rightarrow +\infty}{\longrightarrow } 0,\ \forall\ i\in \{1,\ldots ,k\}\ si\ k\neq 1,\nonumber\\
 && R_n^1\underset{n\rightarrow +\infty}{\longrightarrow } 0,\ if\ k=1.
 \end{eqnarray}
4) There exists $\bar{C}>0$ such that
\begin{eqnarray}
\label{regdecprop3.3}
&&\forall B_s (y)\subset M\backslash \bigcup_{i=1}^k B_{R_n^i}(x_n^i),\ if\ \int_{B_s (y)}e^{v_n}dV \geq \dfrac{\pi}{\rho} \nonumber\\
&&then \ s\geq \bar{C} d_n (y),\ where\ d_n(y)=\displaystyle\inf_{i\in \{1,\ldots , k\}}|y-x_n^i|.
\end{eqnarray}

If $(v_n)_n$ is not compact, from Proposition \ref{mp2}, there exists $(x_n^1 ,r_n^1)_n$ a blow-up for $(v_n)_n$. So there exists a sequence $(b_n)_n$ such that $b_n\underset{n\rightarrow +\infty}{\longrightarrow } +\infty$ arbitrarily slowly, $b_n r_n^1\underset{n\rightarrow +\infty}{\longrightarrow } 0$ and
\begin{equation}
\label{regdecprop3.4}
\lim_{n\rightarrow +\infty}\int_{B_{b_n r_n^1}(x_n^1)}e^{v_n}dV= \dfrac{8\pi}{\rho}.
\end{equation}
Since $b_n\underset{n\rightarrow +\infty}{\longrightarrow } +\infty$ arbitrarily slowly, we can choose an other sequence $(\tilde{b}_n)_n$ such that $\tilde{b}_n \underset{n\rightarrow +\infty}{\longrightarrow } +\infty$ arbitrarily slowly, $\tilde{b}_n\leq \dfrac{b_n}{2}$ and
\begin{equation}
\label{regdecprop3.5}
\lim_{n\rightarrow +\infty}\int_{B_{\tilde{b}_n r_n^1}(x_n^1)}e^{v_n}dV= \dfrac{8\pi}{\rho}.
\end{equation} 
Thanks to (\ref{regdecprop3.4}) and (\ref{regdecprop3.5}), we have
\begin{equation}
\label{24maie10}
\lim_{n\rightarrow +\infty}\int_{B_{b_n r_n^1}(x_n^1) \backslash B_{\tilde{b}_n r_n^1}(x_n^1)}e^{v_n}dV= 0.
\end{equation}
Now, setting  $2R_n^1=b_n r_n^1$, we have, using (\ref{regdecprop3.4}), \eqref{24maie10} and since $\tilde{b}_n \leq \dfrac{b_n }{2}$,
\begin{equation}
\label{regdecprop3.6}
\lim_{n\rightarrow +\infty}\int_{B_{2R_n^1}(x_n^1)\backslash B_{R_n^1}(x_n^1)} e^{v_n}dV\leq \lim_{n\rightarrow +\infty}\int_{B_{b_n r_n^1}(x_n^1) \backslash B_{\tilde{b}_n r_n^1}(x_n^1)}e^{v_n}dV= 0.
\end{equation}
Let $k$ be the biggest integer such that there exist $k$ blow-up $(x_n^1 ,r_n^1)_n,\ldots,(x_n^k ,r_n^k)_n$ and $k$ sequences $(R_n^1)_n,\ldots,(R_n^k)_n>0$ satisfying \eqref{regdecprop2.1.128mai}, (\ref{regdecprop3.1}) and (\ref{regdecprop3.2}). From (\ref{regdecprop3.2}), we have $B_{2R_n^i}(x_n^i)\cap B_{2R_n^j}(x_n^j)=\varnothing$, for all $i\neq j\in \{1,\ldots ,k\}$. We deduce, using (\ref{regdecprop3.1}), that
\begin{equation}
\label{regdecprop3.7}
\lim_{n\rightarrow +\infty}\int_{\bigcup_{i=1}^k B_{2R_n^i}(x_n^i)}e^{v_n}dV= \dfrac{8\pi k}{\rho}.
\end{equation}
On the other hand, one has
\begin{equation}
\label{regdecprop3.8}
\int_{\bigcup_{i=1}^k B_{2R_n^i}(x_n^i)}e^{v_n}dV \leq \int_M e^{v_n}dV=1.
\end{equation}
It follows from (\ref{regdecprop3.7}) and (\ref{regdecprop3.8}) that $k\leq [\dfrac{\rho}{8\pi}]$.\newline
Let's prove that (\ref{regdecprop3.3}) holds. Suppose, by contradiction, that (\ref{regdecprop3.3}) doesn't hold then there exist a sequence of points $(y_n)_n$ and a sequence of positive real numbers $(s_n)_n$ such that $B_{s_n}(y_n)\in M\backslash \displaystyle\bigcup_{i=1}^k B_{R_n^i}(x_n^i)$,
$$\dfrac{s_n}{|y_n-x_n^i|}\underset{n\rightarrow +\infty}{\longrightarrow } 0,\ \forall i\in \{1,\ldots ,k\}$$
and
$$\int_{B_{s_n}(y_n)}e^{v_n}dV \geq \dfrac{\pi}{\rho}.$$ 
Hence, there exist a sequence of points $(x_n)_n$ and a sequence of positive real numbers $(r_n)_n$ such that $B_{r_n}(x_n)\subset M\backslash \displaystyle\bigcup_{i=1}^k B_{R_n^i}(x_n^i)$ and
\begin{equation}
\label{regdecprop3.star}
\int_{B_{r_n}(x_n)}e^{v_n}dV =\max_{B_{r_n}(y)\subset M\backslash \bigcup_{i=1}^k B_{R_n^i}(x_n^i)}\int_{B_{r_n}(y)}e^{v_n}dV=\dfrac{\pi}{\rho}.
\end{equation}
We claim that $\dfrac{r_n}{d_n}\underset{n\rightarrow +\infty}{\longrightarrow } 0$ where $d_n=\displaystyle\inf_{i\in \{1,\ldots , k\}}|x_n -x_n^i|$.\newline
Suppose that $\dfrac{r_n}{d_n}\underset{n\rightarrow +\infty}{\not\longrightarrow } 0$. In this case, we will prove that 
$$\lim_{n\rightarrow +\infty}\int_{B_{r_n}(x_n)}e^{v_n}dV=0.$$
This will lead to a contradiction with \eqref{regdecprop3.star}. Since we assume that $\dfrac{r_n}{d_n}\underset{n\rightarrow +\infty}{\not\longrightarrow } 0$, there exists a constant $C_0>0$ such that
\begin{equation}
\label{regdecprop3.2star}
d_n\geq r_n \geq C_0 d_n.
\end{equation}
We can suppose, up to relabelling the $(x_n^i)_n$, that $d_n=|x_n -x_n^1|$. Let $X_1$ be the set of points such that 
$$X_1= \left\{ (x_n^i)_{n} :\exists C_i >0 \ such\ that\ |x_n^i -x_n^1|\leq C_i d_n \right\}.$$
Up to relabelling the $(x_n^i)_{n}$, $i\neq 1$, we can assume that
$$X_1= \left\{(x_n^1)_{n}, (x_n^{2})_{n}, \ldots , (x_n^{l})_{n}  \right\},$$
with $1 \leq l\leq k$. Let 
\begin{equation}
\label{24maie11}
C=\max\left\{C_i\ :\ i\in \{1,\ldots ,l\}\right\}.
\end{equation}
Thus, we have
$$|x_n^1 -x_n^i|\leq Cd_n, \ \forall i=1,\ldots , l.$$
We will consider two cases depending if $l=1$ or $l>1$.\newline
\textit{First case : $l=1$:}\newline
One can check that the hypothesis of Lemma \ref{regdecprop1} hold for $R_n=R_n^1$ and $S_n=8d_n$. This yields to
$$\lim_{n\rightarrow +\infty}\int_{B_{4d_n}(x_n^1)\backslash B_{2R_n^1}(x_n^1)}e^{v_n}dV= 0.$$
From (\ref{regdecprop3.6}),we get that
\begin{equation}
\label{5juin2012e2}
\lim_{n\rightarrow +\infty}\int_{B_{4d_n}(x_n^1)\backslash B_{R_n^1}(x_n^1)}e^{v_n}dV= 0.
\end{equation}
Since $B_{r_n}(x_n)\subseteq B_{4d_n}(x_n^1)\backslash B_{R_n^1}(x_n^1)$, 
\eqref{5juin2012e2} implies
$$\lim_{n\rightarrow +\infty}\int_{B_{r_n}(x_n)}e^{v_n}dV= 0.$$
This yields to a contradiction with (\ref{regdecprop3.star}). Therefore we deduce that $\dfrac{r_n}{d_n}\underset{n\rightarrow +\infty}{\longrightarrow } 0$.\newline
\textit{Second case : $l >1$.}\newline
Similarly to \eqref{regdecprop3.6}, we see that, for each $(x_n^i)_n$, $i \in \{1,\ldots , l\}$, there exists a sequence $(R_n^i)_n$ such that
$$\dfrac{R_n^i}{r_n^i}\underset{n\rightarrow +\infty}{\longrightarrow } +\infty ,$$
and
\begin{equation}
\label{12nov2012e10}
\lim_{n\rightarrow +\infty}\int_{B_{2R_n^i}(x_n^i)\backslash B_{R_n^i}(x_n^i)} e^{v_n}dV=0.
\end{equation}
We will apply Lemma \ref{regdecprop2} to sequences $(x_n^1)_n,\ldots, (x_n^l)_n$. We notice that, from \eqref{regdecprop3.2}, we have, for all $i\in \{1,\ldots ,l\} $,
$$\dfrac{R_n^i}{\displaystyle\inf_{j\neq i \in \{1,\ldots ,l\}} |x_n^i -x_n^j|}\underset{n\rightarrow +\infty}{\longrightarrow } 0.$$
Now we claim that there exists a constant $\bar{C}>0$ such that
\begin{eqnarray}
\label{25maiee2}
&&\forall B_s (y)\subset \bigcup_{i=1}^l B_{16\tilde{C}d_n}(x_n^i)\backslash \bigcup_{i=1}^l B_{R_n^i}(x_n^i),\ if\ \int_{B_s (y)}e^{v_n}dV \geq \dfrac{\pi}{\rho} \nonumber \\
&&then \ s\geq \bar{C} d_n (y)=\bar{C} \displaystyle\inf_{i\in \{1,\ldots , l\}}|y-x_n^i|,
\end{eqnarray}
where $\tilde{C}=\max \left\{\dfrac{1}{2},C\right\}$, $C$ is the constant defined in \eqref{24maie11}. Let $$B_s (y)\subset \bigcup_{i=1}^l B_{16\tilde{C}d_n}(x_n^i)\backslash \bigcup_{i=1}^l B_{R_n^i}(x_n^i)$$ be such that $\displaystyle\int_{B_s (y)}e^{v_n}dV \geq \dfrac{\pi}{\rho}$. From \eqref{regdecprop3.star}, we deduce that $s\geq r_n\geq C_0 d_n$. On the other hand, since $B_s(y)\subset \displaystyle\bigcup_{i=1}^l B_{16\tilde{C}d_n}(x_n^i)$, we have $\displaystyle\inf_{i\in \{1,\ldots,l\}}|y-x_n^i|\leq 16\tilde{C} d_n$. Thus we get that $s\geq \dfrac{C_0}{16\tilde{C}}\displaystyle\inf_{i\in \{1,\ldots , l\}}|y-x_n^i|$. Therefore, we can take $\bar{C}=\dfrac{C_0}{16 \tilde{C}}$. So \eqref{25maiee2} holds.\newline
Set $\tilde{D}_n= \underset{i,j\in \{1,\ldots ,l\}}{\max}|x_n^i -x_n^j|$. We note that, since $C$ is defined in \eqref{24maie11}, we have
\begin{equation}
\label{25maiee1}
Cd_n \leq \tilde{D}_n \leq 2Cd_n.
\end{equation}
This implies that $4\tilde{D}_n \leq 8Cd_n$. Thus we can apply Lemma \ref{regdecprop2}. This yields to
\begin{equation}
\label{6avrillem1012e1mai1}
\lim_{n\rightarrow +\infty}\int_{\cup_{i=1}^l B_{2\tilde{D}_n}(x_n^i)\backslash \cup_{i=1}^l B_{2R_n^i}(x_n^i)}e^{v_n}dV= 0.
\end{equation}
From \eqref{25maiee2} and since, from \eqref{25maiee1}, $8\dfrac{\tilde{C}}{C}\tilde{D}_n \leq 16\tilde{C}d_n$, reasoning in the same way as for obtaining \eqref{avril2012avantderlem1} of Lemma \ref{regdecprop2}, we have
\begin{equation}
\label{6avrillem1012e1mai2}
\lim_{n\rightarrow +\infty}\int_{\cup_{i=1}^l B_{4\frac{\tilde{C}}{C}\tilde{D}_n}(x_n^i)\backslash \cup_{i=1}^l B_{2\tilde{D}_n}(x_n^i)}e^{v_n}dV= 0.
\end{equation}
It follows, from \eqref{6avrillem1012e1mai1} and \eqref{6avrillem1012e1mai2}, that
\begin{equation}
\label{6avrillem1012e1mai}
\lim_{n\rightarrow +\infty}\int_{\cup_{i=1}^l B_{4\frac{\tilde{C}}{C}\tilde{D}_n}(x_n^i)\backslash \cup_{i=1}^l B_{2R_n^i}(x_n^i)}e^{v_n}dV= 0.
\end{equation}
From \eqref{25maiee1}, we deduce that 

\begin{equation}
\label{6avrillem1012e1}
\lim_{n\rightarrow +\infty}\int_{\cup_{i=1}^l B_{4\tilde{C}d_n}(x_n^i)\backslash \cup_{i=1}^l B_{2R_n^i}(x_n^i)}e^{v_n}dV= 0.
\end{equation}
Since, by \eqref{regdecprop3.2star}, $r_n \leq d_n$, we have $|x_n-x_n^1|+r_n \leq 2d_n \leq 4\tilde{C}d_n$. Thus we deduce that
\begin{equation}
\label{6avrillem1012e2}
\displaystyle B_{r_n}(x_n)\subset \displaystyle \bigcup_{i=1}^l B_{4\tilde{C}d_n}(x_n^i)\backslash \displaystyle\bigcup_{i=1}^l B_{R_n^i}(x_n^i).
\end{equation}
On the other hand, from (\ref{12nov2012e10}), we have, for all $i\in \{1,\ldots ,l\}$,
\begin{equation}
\label{6avrillem1012e3}
\lim_{n\rightarrow +\infty}\int_{B_{2R_n^i}(x_n^i)\backslash B_{R_n^i}(x_n^i)}e^{v_n}dV = 0.
\end{equation}
Combining \eqref{6avrillem1012e1}, \eqref{6avrillem1012e2} and \eqref{6avrillem1012e3}, we deduce that

$$\lim_{n\rightarrow +\infty}\int_{B_{r_n}(x_n)}e^{v_n}dV= 0.$$
This yields to a contradiction with (\ref{regdecprop3.star}). Therefore, we have proved that
\begin{equation}
\label{6avrillem1012e4}
\dfrac{r_n}{d_n}\underset{n\rightarrow +\infty}{\longrightarrow } 0.
\end{equation}
Now, it is easy to see, using \eqref{regdecprop3.star} and proceeding as in the proof of Proposition \ref{mp2}, that $(x_n^{k+1})_n=(x_n)_n$ and $(r_n^{k+1})_n=(r_n)_n$ are a blow-up for $(v_n)_n$. Moreover, since $(x_n ,r_n)_n$ is a blow-up, there exists a sequence $(b_n)_n$, $b_n\underset{n\rightarrow +\infty}{\longrightarrow } +\infty$ arbitrarily slowly such that
$$\lim_{n\rightarrow +\infty}\int_{B_{b_nr_n}(x_n)}e^{v_n}dV=\dfrac{8\pi}{\rho}.$$
So, by setting $2R_n^{k+1}=b_n r_n^{k+1}$, we have
$$\lim_{n\rightarrow \infty}\int_{B_{2R_n^{k+1}}(x_n^{k+1})}e^{v_n}dV = \dfrac{8\pi}{\rho}.$$
From \eqref{6avrillem1012e4} and since $(b_n)_n$, $b_n\underset{n\rightarrow +\infty}{\longrightarrow } +\infty$ arbitrarily slowly, we also can assume without loss of generality that
$$\dfrac{R_n^{k+1}}{\displaystyle\inf_{i\in \{1, \ldots , k\}}|x_n^{k+1}-x_n^i|}\underset{n\rightarrow +\infty}{\longrightarrow } 0.$$
Therefore, the sequences $(x_n^1,r_n^1)_n,\ldots, (x_n^{k+1},r_n^{k+1})_n$ and $(R_n^1)_n,\ldots ,(R_n^{k+1})_n$ satisfy \eqref{regdecprop2.1.128mai}, (\ref{regdecprop3.1}) and (\ref{regdecprop3.2}). We obtain, this way, a contradiction with the maximality of $k$. Therefore, there exist $k$ blow-up $$(x_n^1, r_n^1)_{n},\ldots ,(x_n^k, r_n^k)_{n}$$ with $1\leq k\leq [\dfrac{\rho}{8\pi}]$ (where $[\dfrac{\rho}{8\pi}]$ stands for the integer part of $\dfrac{\rho}{8\pi}$), and there exist $k$ sequences $(R_n^1)_n,\ldots ,(R_n^k)_n >0$ satisfying \eqref{regdecprop2.1.128mai}, (\ref{regdecprop3.1}) and (\ref{regdecprop3.2}). If $k=1$, then, thanks to \eqref{regdecprop2.1.128mai}, (\ref{regdecprop3.1}) and (\ref{regdecprop3.2}), the hypothesis of Lemma \ref{regdecprop1} for $x_n =x_n^1$, $R_n=R_n^1$ and $S>0$ a real number such that $M\subset B_{\frac{S}{2}}(x_n^1)$ hold. Thus, we have
$$\lim_{n\rightarrow \infty} \int_{M\backslash B_{2R_n^1}(x_n^1)} e^{v_n}dV=0.$$
This establishes Theorem \ref{regdecthm3} in the case $k=1$.\\
If $k>1$, then, thanks to \eqref{regdecprop2.1.128mai}, (\ref{regdecprop3.1}) and (\ref{regdecprop3.2}), the hypothesis of Lemma \ref{regdecprop2ter} hold. Hence, we obtain that
$$\displaystyle\lim_{n\rightarrow +\infty} \int_{ M\backslash \bigcup_{i=1}^k B_{2R_n^i}(x_n^i)}e^{v_n}dV=0.$$
This concludes the proof of Theorem \ref{regdecthm3}.
\section{Convergence of the flow.}
This section is devoted to the proof of Theorems \ref{conv1}, \ref{conv2} and \ref{conv3}.
\subsection{Proof of Theorem \ref{conv1}.}
In this subsection, we will prove Theorem \ref{conv1}. Let $v : M\times [0,+\infty )\rightarrow \R$ be the global solution of \eqref{E:flot}. In all this subsection, we will assume without loss of generality that $\displaystyle\int_M e^{v(t)}dV=1$, $\forall t\geq 0$. In the following, $C$ will denote constants not depending on $t$.\\

In order to prove Theorem \ref{conv1}, we need to bound uniformly $\left\|v(t)\right\|_{H^2(M)}$, $t\geq 0$, in time. For this, we will first bound $\left\|v(t)\right\|_{H^1(M)}$, $t\geq 0$, uniformly in time. To bound $\left\|v(t)\right\|_{H^1(M)}$, we will use the compactness Corollary \ref{regdecthm2}. More precisely, by using Corollary \ref{regdecthm2}, we will first prove that there exists a sequence $(t_n)_n$, $t_n\underset{n\rightarrow +\infty}{\longrightarrow}+\infty$, such that
$$\left\|v(t_n)\right\|_{H^2 (M)}\leq C, \ \forall n\geq 0.$$
Therefore we aim to prove that there exists a sequence $(t_n)_n$, $t_n\underset{n\rightarrow +\infty}{\longrightarrow}+\infty$, such that, setting $v_n=v(t_n)$ and $h_n=-\dfrac{\partial v}{\partial t}(t_n)$, the sequence $(v_n)_n\subseteq H^2(M)$ satisfies conditions \eqref{regdeceq1} of Corollary \ref{regdecthm2}. First, we show that there exists a sequence $(t_n)_n$, $t_n\underset{n\rightarrow +\infty}{\longrightarrow}+\infty$, such that \eqref{regdeceq1}(ii) is satisfied for $v_n=v(t_n)$. We recall that, for all $T>0$, 
$$\int_{0}^{T}\int_M\left(\dfrac{\partial v(t)}{\partial t}\right)^{2}e^{v(t)} dVdt=J_\rho(v(0))-J_\rho(v(T)).$$
Using hypothesis \eqref{19juin2012e1}, we deduce that there exists a sequence $(t_n)_n$ such that $n\leq t_n \leq n+1$, for all $n\in \N$, and
\begin{equation}
\label{propsup1.1}
\underset{n\rightarrow +\infty}{\lim}\int_M |\dfrac{\partial v(t_n)}{\partial t}|^2e^{v(t_n)}dV= 0.
\end{equation}
The following proposition shows that conditions \eqref{regdeceq1}(i) of Corollary \ref{regdecthm2} is satisfied.
\begin{prop}
We have
\begin{equation}
\label{mp4}
-\dfrac{\partial e^{v(x,t)}}{\partial t}+\rho e^{v(x,t)}\geq -C,\ \forall t\geq 0\ and\ x\in M.
\end{equation}
\end{prop}
\begin{proof}
We set
 $$ R(x,t)=e^{-v(x,t)}\left(- \Delta v(x,t)+Q(x)\right).$$
We can rewrite the equation \eqref{E:flot} satisfied by $v$ in the following way
 $$\dfrac{\partial v(x,t)}{\partial t}=-\left(R(x,t)- \rho \right). $$
Hence
 \begin{eqnarray*}
 \dfrac{\partial R(x,t)}{\partial t}  & =& R(x,t) \left(R(x,t) - \rho \right)+e^{-v(x,t)}\Delta R(x,t). 
 \end{eqnarray*}
 Set $R_{\min}\left(t\right)=\displaystyle\min_{x \in M}R\left(x,t\right)$. Using the maximum principle, we find 
 
 \begin{eqnarray*}
 \dfrac{\partial R_{\min}(t)}{\partial t} & \geq &- \rho   R_{\min}(t) .
 \end{eqnarray*}
Integrating between $0$ and $t$, we have
 $$R_{\min}(t)\geq \exp\left({-\rho t}\right)  R_{\min}(0) .$$
This implies that
\begin{eqnarray}
\label{mei}
&&-\dfrac{\partial e^{v(x,t)}}{\partial t}+ \rho e^{v(x,t)}\nonumber \\
&\geq &-\left|R_{\min}(0) \right|\exp\left({-\rho t}+v(x,t)\right).
\end{eqnarray}
Set $v_{\max}(t)=\displaystyle\max_{x \in M}v\left(x,t\right)$. By the maximum principle, we have
 
 \begin{eqnarray*}
 \dfrac{\partial }{\partial t}e^{v_{\max}(t)}  & \leq &\rho\left(\frac{1}{\rho} \left\|Q\right\|_{L^\infty (M)} +e^{v_{\max}(t)}  \right).
 \end{eqnarray*}
Integrating between $0$ and $t$, we get
  \begin{equation}
 \label{mp4.1}
 e^{v_{\max}\left(t\right)-\rho t}\leq e^{v_{\max}\left(0\right)}+\frac{1}{\rho}\left\|Q\right\|_{L^\infty (M)}  -\frac{1}{\rho}\left\|Q\right\|_{L^\infty (M)}e^{-\rho t}\leq C.
 \end{equation}
Combining (\ref{mei}) and (\ref{mp4.1}), we finally obtain
 $$-\dfrac{\partial e^{v(x,t)}}{\partial t}+\rho e^{v(x,t)}\geq -C \left| R_{\min}\left(0\right)\right|\geq -C.$$
\end{proof}

We are now in position to bound $\left\|v(t)\right\|_{H^1(M)}$, $t\geq 0$, uniformly in time.
\begin{prop}
Let $\rho \in (8k\pi, 8(k+1)\pi )$, $k\in \N^\ast$ and $v(t):M\rightarrow \R$ be the solution of \eqref{E:flot}. Suppose that
\begin{equation}
\label{19juin2012e1}
J_\rho (v(t))\geq -C,\ \forall t\geq 0.
\end{equation}
Then we have
\begin{equation}
\label{11nov2012e3}
\left\|v(t)\right\|_{H^1(M)}\leq C,\ \forall t \geq 0.
\end{equation}
\end{prop}
\begin{proof}
\setcounter{etape}{0}
 
Thanks to \eqref{propsup1.1} and \eqref{mp4}, by using Corollary \ref{regdecthm2}, there exists a constant $C>0$ such that
$$ \left\| v(t_n)\right\|_{H^2 (M)}\leq C,$$
where $(t_n)_n$ is the sequence defined in \eqref{propsup1.1}. By Sobolev's embedding Theorem, it follows that $\left\|v(t_n)\right\|_{C^\alpha (M)}\leq C$, $\forall \alpha\in (0,1)$. Since $t_n\underset{n\rightarrow +\infty}{\longrightarrow} +\infty$, for all $t\geq 0$ large enough, there exists $n\in \N$ such that $t_n \leq t \leq t_{n+1}$. Moreover, since $|t_{n+1}-t_n|\leq 2$, we have $|t-t_n|\leq 2$. We claim that, for all $p>1$,
\begin{equation}
\label{19juin2012e4}
\int_M e^{pv(t)}dV \leq C,\ \forall t\geq 0.
\end{equation}
Since $v(t)$ satisfies \eqref{E:flot}, integrating by parts and using Young's inequality, we see that
\begin{eqnarray*}
\dfrac{\partial}{\partial t}\int_M e^{pv(t)}dV &=& -p(p-1)\int_M |\nabla v(t)|^2  e^{(p-1)v(t)}dV-p \int_M Q e^{(p-1)v(t)}dV\\
&+& p \frac{\rho}{a} \int_M e^{pv(t)}dV \\
& \leq & C \int_M e^{(p-1)v(t)}dV+p \frac{\rho}{a} \int_M e^{pv(t)}dV \\
&\leq & C +C \int_M e^{pv(t)}dV.
\end{eqnarray*}
Setting $y(t)=\displaystyle\int_M e^{pv(t)}dV$ and integrating the previous inequality between $t_n$ and $t$, it follows that
\begin{equation*}
y(t) \leq e^{C(t-t_n)}y(t_n)+C\left(e^{C(t-t_n)}-1 \right).
\end{equation*}
Since $\left\|v(t_n)\right\|_{C^\alpha (M)}\leq C$, $\alpha\in (0,1)$, and $|t-t_n|\leq 2$, we have that \eqref{19juin2012e4} is satisfied. Let's fix $t\geq 0$ and set  $$M_\varepsilon = \left\{x \in M : e^{v(x,t)}<\varepsilon\right\},$$

\noindent where $\varepsilon > 0$ is a real number which will be determined later. We have
\begin{eqnarray}
\label{2ee21.1}
1= \int_{M}e^{v(t)}dV & =& \int_{M_{\varepsilon}}e^{v(t)}dV+\int_{M\backslash M_{\varepsilon}} e^{v(t)}dV \nonumber\\
& \leq & \varepsilon |M_{\varepsilon}|+ \left|M\backslash M_{\varepsilon}\right|^{1-\frac{1}{p}} \left(\int_{M}e^{pv(t)}dV\right)^{\frac{1}{p}}.
\end{eqnarray}
Taking $\varepsilon=\frac{1}{2|M|}$ and from \eqref{19juin2012e4}, we deduce that
$$\frac{1}{2}\leq C |M\backslash M_{\varepsilon}|^{1-\frac{1}{p}}.$$
Since $p>1$, we get that
\begin{equation}
\label{3juillet2012e1}
|M\backslash M_{\varepsilon}|\geq \left(\dfrac{1}{2C}\right)^{\frac{p}{p-1}}>0.
\end{equation}
We set $A=M\backslash M_\varepsilon$. By definition of $A$, we have 
\begin{equation}
\label{rmqe1}
\int_{A}v(t)dV \geq \ln \left(\dfrac{1}{2|M|}\right) |A|.
\end{equation}
On the other hand, we have 
$$\int_{A}v(t) dV \leq \int_{A}e^{v(t)} dV \leq 1.$$ 
We deduce from the previous inequality and (\ref{rmqe1}) that there exists a constant $C$ such that
\begin{equation}
\label{rmqe2}
\left|\int_{A}v(t)dV  \right|\leq C.
\end{equation} 
From \eqref{3juillet2012e1} and (\ref{rmqe2}), arguing the same way as in Proposition \ref{etape4}, we deduce that, for all $t\geq 0$,
$$\left\|v(t)\right\|_{H^1 (M)}\leq C.$$
\end{proof}

\vspace{12pt}
\vspace{12pt}

\noindent\textit{Proof of Theorem \ref{conv1}.} First, let us prove that
\begin{equation}
\label{11nov2012e4}
\int_M \left(\Delta v(t)\right)^2 dV \leq C,\ \forall t\geq 0.
\end{equation}
We follow Brendle \cite{MR1999924} arguments. We set $$V(t)=\dfrac{\partial v(t)}{\partial t}$$ and $$y(t)=\int_{M}V^{2}(t)e^{v(t)}dV.$$
We claim that $y(t)\underset{t\rightarrow +\infty}{\longrightarrow} 0$. By \eqref{11nov2012e3}, we have, for all $T\geq 0$,
\begin{equation}
\label{conve.1}
\int_0^T\int_M \left(\dfrac{\partial v(t)}{\partial t}\right)^2 e^{v(t)} dVdt \leq J_\rho (v(0))-J_\rho (v(T))\leq C,
\end{equation}
where $C$ is a constant not depending on $T$. Let $\varepsilon$ be some positive real number. From \eqref{conve.1}, we deduce that there exists $t_{0}\geq 0$ such that $y(t_{0})\leq \varepsilon$.\newline
We want to prove that
$$y(t)\leq 3\varepsilon , \ \forall t\geq t_{0}.$$
Otherwise, we define
$$t_{1}=\inf \left\{\ t\geq t_{0}\ :\ y(t)\geq 3\varepsilon \right\}<+\infty.$$
This implies that
\begin{equation}
\label{11nov2012e5}
y(t)\leq 3\varepsilon, \ \forall t_{0}\leq t\leq t_{1}.
\end{equation}
Since $\dfrac{\partial v(t)}{\partial t}=e^{-v(t)}\left(\Delta v(t) -Q\right)+\rho $, we arrive, by using \eqref{11nov2012e5}, at 
\begin{equation}
\label{ie3.1}
\int_{M}e^{-v(t)}\left(\Delta v(t) -Q\right)^{2}dV=y(t)+\rho^2\leq C_1,\ \forall t_{0}\leq t\leq t_{1},
\end{equation}
where, in the following, $C_1$ will denote constant depending on $\varepsilon$ and thus on $t_1$. From \eqref{19juin2012e4}, we have, for all $t\geq 0$,
\begin{equation}
\label{ie3.2}
\int_{M}e^{3v(t)}dV\leq C,
\end{equation}
where, in the following, $C$ will denote constant not depending on $t_1$. Using H\"{o}lder's inequality, (\ref{ie3.1}) and (\ref{ie3.2}), we obtain, $\forall t_{0}\leq t\leq t_{1}$,
\begin{eqnarray*}
\int_{M}\left|\Delta v(t) -Q\right|^{\frac{3}{2}}dV & \leq &\left(\int_{M}e^{-v(t)}\left(\Delta v(t) -Q\right)^{2}dV  \right)^{\frac{3}{4}}\left(\int_{M}e^{3v(t)}dV \right)^{\frac{1}{4}} \\
& \leq & C_1, 
\end{eqnarray*}
thus, $\forall t_{0}\leq t\leq t_{1}$, one has $\displaystyle \int_{M}\left|\Delta v(t)\right|^{\frac{3}{2}}dV\leq C_1.$
From Sobolev's embedding Theorem, we get that
\begin{equation}
\label{3juillet2012e2}
\left|v(t)\right|\leq C_1,\ \forall t_{0}\leq t\leq t_{1}.
\end{equation}
On the other hand, we see that $V(t)=\dfrac{\partial v(t)}{\partial t}$ satisfies
\begin{eqnarray}
\label{19juin2012e10}
\frac{\partial V(t)}{\partial t}&=& -V(t) e^{-v(t)}\Delta v(t)+ e^{-v(t)}\Delta V(t) +QV(t)e^{-v(t)}.
\end{eqnarray}
Now using \eqref{19juin2012e10}, we have, for all $t_{0}\leq t\leq t_{1}$,
\begin{eqnarray*}
\frac{\partial y(t)}{\partial t} & =& \frac{\partial }{\partial t}\left(\int_{M}V^{2}(t)e^{v(t)}dV\right) \nonumber\\
& =& 2\int_{M}V(t)e^{v(t)}\left(e^{-v(t)}\Delta V(t)- V(t) e^{-v(t)} \Delta v(t) +QV(t)e^{-v(t)} \right) dV\nonumber \\
&+&\int_{M}V^{3}(t)e^{v(t)}dV \nonumber
\end{eqnarray*}
Integrating by parts, we obtain
\begin{equation}
\label{conve.2}
\frac{\partial y(t)}{\partial t} =-2 \int_{M}\left|\nabla V(t)\right|^{2}dV- \int_{M}V^{3}(t)e^{v(t)}dV+2\rho\int_{M}V^{2}(t)e^{v(t)}dV. 
\end{equation}
Using Gagliardo-Nirenberg's inequality, we get
$$\left\|V(t)\right\|_{L^{3}_{g_{1}}\left(M\right)}\leq C \left\|V(t)\right\|_{L^{2}_{g_{1}}\left( M\right)}^{\frac{2}{3}}\left\|V(t)\right\|_{H^{1}_{g_1}\left(M\right)}^{\frac{1}{3}},$$
where the norms are taken with respect to the metric $g_1(t)=e^{v(t)}g$. From \eqref{3juillet2012e2}, we notice that the first eigenvalue of the laplacian $\tilde{\lambda}_1(t)$ with respect to the metric $g_1(t)$ satisfies, $\forall t_0\leq t \leq t_1$,
\begin{equation}
\label{raleygh}
\tilde{\lambda}_1 (t)\geq C_1.
\end{equation}
From the fact that $\displaystyle\int_{M}Ve^{v}dV=0$, Poincare's inequality and \eqref{raleygh}, we have
\begin{equation}
\label{conve.3}
\int_{M}e^{v}\left|V\right|^{3}dV\leq C_1 \left(\int_{M}V^{2}e^{v}dV\right)\left(\int_{M}\left|\nabla V\right|^{2}dV\right)^{\frac{1}{2}}.
\end{equation}
Thus we obtain, from (\ref{conve.2}), (\ref{conve.3}) and Young's inequality,
$$\frac{\partial}{\partial t}\left(\int_{M}V^{2}e^{v}dV\right)\leq C_1 \left(\int_{M}V^{2}e^{v}dV\right)^{2}+C\left(\int_{M}V^{2}e^{v}dV\right),$$
i.e.
$$\frac{\partial}{\partial t}y(t)\leq C_1 y^{2}(t)+Cy(t).$$
Since $y(t_0) \leq \varepsilon$ and $y(t_1)=3\varepsilon$, we find
$$2\varepsilon \leq y(t_{1})-y(t_{0})\leq \left(C_1+C\right) \int_{t_{0}}^{t_{1}}y(t)dt.$$
Choosing $t_{0}$ large enough, we have $\left(C_1+C\right)\displaystyle\int_{t_{0}}^{+\infty}y(t)dt\leq \varepsilon , $
and thus we obtain a contradiction. We conclude that
$$y(t)\underset{t\rightarrow +\infty}{\longrightarrow} 0.$$
Thereby, $t_1=+\infty$. This implies that all estimates we previously got, hold for all $t\geq 0$. Thus, we have, for all $t\geq 0$, $|v(t)|\leq C,$ and
$$\int_M e^{-v(t)}(\Delta v(t) - Q)^2 dV \leq C.$$
It follows that, for all $t\geq 0$, $\displaystyle\int_M \left(\Delta v(t)\right)^2 dV \leq C$. Thus, using \eqref{11nov2012e3}, we have, for all $t\geq 0$, $\left\|v(t) \right\|_{H^2 (M)}\leq C.$ Therefore, there exist a function $v_\infty \in H^{2}(M)$ and a sequence $(t_n)_n$, $t_n\underset{n\rightarrow +\infty}{\longrightarrow}+\infty$ such that
$$v(t_n) \underset{n\rightarrow +\infty}{\longrightarrow}v_\infty \ weakly\ in\ H^2(M),$$
and
$$v(t_n) \underset{n\rightarrow +\infty}{\longrightarrow}v_\infty \ in\ C^\alpha (M), \ \alpha\in (0,1). $$ 
It is easy to check that $v_\infty$ is a solution to
$$-\Delta v_\infty +Q=\rho\dfrac{e^{v_\infty}}{\int_M e^{v_\infty} dV},$$
and, by boothstrap regularity arguments, we have $v_\infty \in C^\infty (M)$.
To obtain that $\left\|v(t_n)-v_\infty\right\|_{H^2(M)} \underset{n\rightarrow +\infty}{\longrightarrow}0$, we just notice that
\begin{eqnarray*}
\int_M \left(\Delta v(t_n) -\Delta v_\infty \right)^2 dV&=& \int_M \left(\dfrac{\rho}{a}(e^{v_\infty}-e^{v(t_n)})+\dfrac{\partial e^{v(t_n)}}{\partial t}\right)^2dV\\
&\leq & C\int_M \left(e^{v_\infty}-e^{v(t_n)}\right)^2dV+ C \int_M \left|\dfrac{\partial v}{\partial t}(t_n)\right|^2 e^{v(t_n)} dV\\
& \underset{n\rightarrow +\infty}{\longrightarrow} & 0.
\end{eqnarray*} 
Since the flow is a gradient flow for the functional $J_\rho$ which is real analytic, from a general result of Simon \cite{MR727703}, we finally obtain that
$$ \left\|v(t)-v_\infty \right\|_{H^2(M)}\underset{n\rightarrow +\infty}{\longrightarrow} 0. $$

\subsection{Proof of Theorem \ref{conv2}.}
We will prove the existence of an initial data $v_0\in C^\infty (M)$ for the flow \eqref{E:flot} such that the functional $J_\rho (v(t))$, $t\geq 0$, is uniformly bounded from below. From standard parabolic theory, it is easy to see that if $v_0 \in C^\infty (M)$ then the solution $v$ of \eqref{E:flot} belongs to $C^{\infty}(M\times [0,+\infty ))$. Let $X$ be the space of functions $C^\infty(M)$ endowed with the $\left\|.\right\|_{C^{2+\alpha}(M)}$ norm and
 $$\Phi :\left\{\begin{aligned}
 &X&\times &[0,&+\infty )&\longrightarrow &C^{\infty}(M\times [0,+\infty)\\
 &(v&,&t&)&\longrightarrow &\Phi(v,t)
 \end{aligned}\right.$$ 
where $\Phi (v,t)$ is a solution of
$$\left\{\begin{array}{l }
  \dfrac{\partial \Phi}{\partial t}(v,t)=e^{-\Phi (v,t)}\Delta \Phi (v,t)-e^{-\Phi (v,t)}Q+\frac{\rho}{\int_M e^{\Phi (v,t)}dV} \\ \Phi(v,0)=v.
         
\end{array}  \right. $$

\noindent Suppose that, $\forall v \in X$, we have
\begin{equation} \label {Ctrd} J_\rho \left(\Phi (v,t)\right) \underset{t \rightarrow + \infty}{\longrightarrow} -\infty . \end{equation}
Let $L>0$. Using \eqref{Ctrd}, we will prove that $\left\{v\in X : J_\rho (v) \leq -L\right\}$ is contractible. But following the same arguments as in Malchiodi \cite{MR2483132}, one can prove that there exists $L_1>0$ such that $\left\{v\in X : J_\rho (v) \leq -L_1\right\}$ is not contractible. Let us prove that if \eqref{Ctrd} is satisfied then $\left\{v\in X : J_\rho (v) \leq -L\right\}$ is contractible. We proceed in two steps.
\setcounter{etape}{0}
\begin{etape}
\label{etapeborne1}
Let $L>0$ be fixed and  $$ T_v= \\inf\left\{ t\geqslant 0,   J_\rho\left(\Phi (v,t)\right)\leq -L\right\},$$ then the function $T:\left\{\begin{array}{llccc }
  C^{2+\alpha}(M) &\longrightarrow &\R\\  v &\longrightarrow &T_v
         \end{array} \right.$
is continuous.
\end{etape}
\noindent\textit{Proof of Step 1.} From (\ref{Ctrd}), we have
$$\left\{ t\geqslant 0,   J_\rho\left(\Phi (v,t)\right)\leq -L\right\}\neq \varnothing ,$$
and, from the uniqueness of solutions of \eqref{E:flot} having same initial data, one can prove that $J_\rho\left(\Phi (v,t)\right)$ is strictly decreasing on $[0,+\infty)$. Let $\bar{v}\in C^{\infty}(M)$ and $(v_n)_{n}\in C^{\infty}(M)$ be a sequence such that $v_n\underset{n\rightarrow +\infty}{\longrightarrow} \bar{v}$ in $C^{2+\alpha}(M)$, we claim that $T_{v_n}\underset{n\rightarrow +\infty}{\longrightarrow} T_{\bar{v}}$. To prove this, we will consider two cases depending on the value of $J_\rho (\bar{v})$.\newline
\textit{First case.} Suppose that $J_\rho (\bar{v})<-L$. Since the function $t\rightarrow J_\rho(\Phi (\bar{v},t))$ is decreasing, we have $J_\rho (\Phi (\bar{v},t))<-L$ for all $t\geq 0$. We deduce that $T_{\bar{v}}=0$. Since $v_n\underset{n\rightarrow +\infty}{\longrightarrow} \bar{v}$ in $C^{2+\alpha}(M)$, it is easy to see that
$$J_\rho (v_n)\underset{n\rightarrow +\infty}{\longrightarrow} J_\rho (\bar{v}).$$
Thus, there exists $n_0 \in \N$ such that $J_\rho (v_n)\leq -L $ for all $n\geq n_0$. So we obtain that $T_{v_n}=0=T_{\bar{v}}$ for all $n\geq n_0$. This implies that
$$T_{v_n}\underset{n\rightarrow +\infty}{\longrightarrow} T_{\bar{v}}.$$
\textit{Second case.} Suppose that $J_\rho (\bar{v})\geq-L$. In this case, $T_{\bar{v}}$ verifies $J_\rho (\Phi (\bar{v},T_{\bar{v}}))=-L$. Setting $T_n:=T_{v_n}$ and supposing that 
$T_n$ does not converge to $T_{\bar{v}}$, then, up to extracting a sub-sequence, there exists $\varepsilon_0>0$ such that $|T_n-T_{\bar{v}}|\geq \varepsilon_0$. So we have $T_n \geq \varepsilon_0 +T_{\bar{v}}$ or $T_n \leq -\varepsilon_0 +T_{\bar{v}}$. Suppose, without loss of generality, that 
\begin{equation}
\label{jrhoborne}
T_n \geq \varepsilon_0 +T_{\bar{v}}.
\end{equation}
Set $T=T_{\bar{v}}+\varepsilon_0+1$. Since $v_n \underset{n\rightarrow +\infty}{\longrightarrow} \bar{v}$ in $C^{2+\alpha}(M)$, by Proposition \ref{propdonin}, it is easy to see that 
\begin{equation}
\label{jrhoborne2}
J_\rho (\Phi (v_n, t))\underset{n\rightarrow +\infty}{\longrightarrow} J_\rho (\Phi (\bar{v}, t)),
\end{equation}
for all $t$ fixed in $[0,T]$. Since $t\rightarrow J_\rho (\Phi (\bar{v},t))$ is strictly decreasing, we have $$\alpha_1= J_\rho (\Phi (\bar{v}, T_{\bar{v}}))-J_\rho (\Phi (\bar{v}, T_{\bar{v}}+\varepsilon_0))>0.$$

\noindent From (\ref{jrhoborne2}) , we get, since $T_{\bar{v}}+\varepsilon_0 \in [0,T]$,
$$J_\rho (\Phi (v_n, T_{\bar{v}}+\varepsilon_0))\underset{n\rightarrow +\infty}{\longrightarrow} J_\rho (\Phi (\bar{v}, T_{\bar{v}}+\varepsilon_0))=-L-\alpha_1$$
and from (\ref{jrhoborne}),
$$J_\rho(\Phi (v_n, T_n))\leq J_\rho(\Phi (v_n, T_{\bar{v}}+\varepsilon_0)),$$
this implies that, if $n$ tends to $+\infty$, $-L\leq -L-\alpha_1$, thus we obtain a contradiction.

\begin{etape}
If \eqref{Ctrd} holds, then the set $\left\{v\in X: J_\rho(v)\leqslant -L\right\}$ is contractible.
\end{etape}
         
\noindent\textit{Proof of Step 2.} We will construct a deformation retract from $\left\{v\in X\right\}$ into $\left\{v\in X: J_\rho(v)\leqslant -L\right\}$. Since $\left\{v\in X\right\}$ is contractible, then $\left\{v\in X: J_\rho(v)\leqslant -L\right\}$ will also be contractible. We denote by $h$ the one-to-one function defined by
$$\begin{array}{llccc }
h(t):&[0,1) &\longrightarrow &[0, +\infty) \\
&t &\longrightarrow &\frac{t}{1-t},
\end{array}$$
and by $\eta(v,t):X\times [0,1]\rightarrow X$ the function defined by 
$$\eta(v,t) =\left\{\begin{array}{r @{ \textrm { if } } l  }
     \Phi\left(v,h(t)\right)& h(t)\leqslant T_v \\
     \Phi\left(v,T_v\right)& h(t)\geqslant T_v.  \\
\end{array}\right. $$

\noindent In a first time, let us prove that $\eta =\Phi \circ \Phi_1 :X\times [0,1)\rightarrow X$ is continuous where $\Phi_1 :X \times [0,1)\rightarrow X\times [0,+\infty)$ is the function defined by
$$\Phi_1(v,t) =\left\{\begin{array}{r @{ \textrm { if } } l  }
     \left(v,h(t)\right)& h(t)\leqslant T_v \\
     \left(v,T_v\right)& h(t)\geqslant T_v.  \\
\end{array}\right. $$
 From Step \ref{etapeborne1}, $\Phi_1 :X \times [0,1)\rightarrow X\times [0,+\infty)$ is a continuous function. Therefore, to prove that $\eta$ is a continuous application from $X\times [0,1)\rightarrow X$, it is sufficient to prove that, for $T>0$ fixed, $\Phi :X\times [0,T]\rightarrow X$ is continuous. Let $(v_n ,t_n)\in C^\infty (M) \times [0,T]$ such that $v_n \underset{n\rightarrow +\infty}{\longrightarrow} v$ in $C^{2+\alpha}(M)$, where $v\in C^\infty (M)$ and $t_n\underset{n\rightarrow +\infty}{\longrightarrow} t\in [0,T]$. Then, we have
\begin{eqnarray}
\label{9avrilret3}
\left\|\Phi (v_n,t_n)-\Phi (v, t)  \right\|_{C^{2+\alpha}(M)} &\leq & \left\|\Phi (v_n,t_n)-\Phi (v_n, t)  \right\|_{C^{2+\alpha}(M)}\nonumber\\
&+&\left\|\Phi (v_n,t)-\Phi (v, t)  \right\|_{C^{2+\alpha}(M)}
\end{eqnarray}
Since $\Phi (v_n,.)\in C^{\infty}(M\times [0,T])$, we get, from Theorem \ref{exglob}, that, for all $t\in [0,T]$, 
$$\left\|\dfrac{\partial \Phi (v_n,t)}{\partial t}  \right\|_{C^{2+\alpha}(M)}\leq C_T,$$ 
where, in the following, $C_T$ denotes a constant not depending on $n$. We deduce that
\begin{eqnarray}
\label{9avrilret1}
&&\left\|\Phi (v_n,t_n)-\Phi (v_n, t)  \right\|_{C^{2+\alpha}(M)}\nonumber \\
&=&\left\|\int_{t_n}^t \dfrac{\partial \Phi }{\partial s} (v_n, s)ds  \right\|_{C^{2+\alpha}(M)}\nonumber \\
&\leq & \left|t_n -t\right| \max_{s\in [t_n ,t]}\left\|\dfrac{\partial \Phi (v_n,s)}{\partial s}  \right\|_{C^{2+\alpha}(M)}\underset{n\rightarrow +\infty}{\longrightarrow} 0.
\end{eqnarray}
On the other hand, using Proposition \ref{propdonin}, we have, for all $t\in [0,T]$,
\begin{equation}
\label{9avrilret2}
\left\|\Phi (v_n,t)-\Phi (v, t)  \right\|_{C^{2+\alpha}(M)} \leq C_T \left\|v_n -v  \right\|_{C^{2+\alpha}(M)}\underset{n\rightarrow +\infty}{\longrightarrow} 0.
\end{equation}
Combining \eqref{9avrilret3} , \eqref{9avrilret1} and \eqref{9avrilret2}, we find that
$$\left\|\Phi (v_n,t_n)-\Phi (v, t)  \right\|_{C^{2+\alpha}(M)}\underset{n\rightarrow +\infty}{\longrightarrow}0.$$  
Thus $\eta$ is continuous from $X\times [0,1)\rightarrow X$. It remains to prove that it is continuous on $X\times [0,1]$. Let $(v_n,t_n)\in C^\infty (M) \times [0,1]$ such that $v_n\underset{n\rightarrow +\infty}{\longrightarrow}\bar{v}$ in $C^{2+\alpha}(M)$, where $\bar{v}\in C^{\infty}(M)$, and $t_n\underset{n\rightarrow +\infty}{\longrightarrow} 1$. From Step \ref{etapeborne1}, we have $T_{v_n}=T_n\underset{n\rightarrow +\infty}{\longrightarrow} T_{\bar{v}}$. Since $T_n$ is finite and $t_n\underset{n\rightarrow +\infty}{\longrightarrow} 1$, it follows that $h(t_n)\underset{n\rightarrow +\infty}{\longrightarrow}+\infty$, so, for $n$ large enough, $h(t_n)\geq T_n$ and thus $\eta(v_n,t_n)=\Phi(v_n, T_n)$. We have, in the same way as \eqref{9avrilret1} and \eqref{9avrilret2}, that
 \begin{eqnarray*}
&& \left\|\eta (v_n,t_n)-\eta (\bar{v},1)\right\|_{C^{2+\alpha}(M)}= \left\|\Phi (v_n,T_n)-\Phi (\bar{v},T_{\bar{v}})\right\|_{C^{2+\alpha}(M)}\\
 &\leq &\left\|\Phi (v_n,T_n) -\Phi (\bar{v},T_n)\right\|_{C^{2+\alpha}(M)}+\left\|\Phi (\bar{v},T_n)-\Phi (\bar{v},T_{\bar{v}})\right\|_{C^{2+\alpha}(M)}\underset{n\rightarrow +\infty}{\longrightarrow} 0.
 \end{eqnarray*}  
Therefore $\eta$ is continuous from $X\times [0,1]\rightarrow X$. Now, it is easy to check that $\eta$ is a deformation retract from $X$ into $\left\{v\in X: J_\rho(v)\leqslant -L\right\}$. Therefore, we deduce that $\{v\in X : J_\rho\leqslant-L\}$ is contractible.

\subsection{Non-convergence of the flow: proof of Theorem \ref{conv3}.}
To prove Theorem \ref{conv3}, it is sufficient to prove that there exists a real number $C>0$ depending on $M,Q$ and $\rho$ such that $\forall v_{0}\in C^{2+\alpha}(M)$ satisfying $J_\rho(v_{0})\leq -C$, then the solution $v(t)$ of flow \eqref{E:flot}, with $v(x,0)=v_0(x)$, for all $x\in M$, satisfies
$$J_\rho (v(t))\underset{t\rightarrow +\infty}{\longrightarrow}-\infty.$$
We recall (see Li \cite{MR1673972}) that there exists a constant $C_0\geq 0$ depending on $M, Q$ and $\rho$ such that, for any solution $w\in C^{2+\alpha}(M)$, $\alpha\in (0,1)$, of
$$-\Delta w +Q=\frac{\rho e^{w}}{\int_M e^{w}dV},$$
then
\begin{equation}
\label{12nov2012e1}
\left\|w \right\|_{C^{2+\alpha}(M)}\leq C_0.
\end{equation}
Since $J_\rho(v(t))$ is decreasing, if $\displaystyle \lim_{t\rightarrow +\infty} J_\rho (v(t))\neq -\infty$ then there exists $L\in \R$ such that $$J_\rho (v(t))\geq L, \ \forall t\in [0,+\infty).$$ 
From Theorem \ref{conv1}, there exists a function $v_\infty \in C^\infty (M)$ such that $$\left\|v(t)-v_\infty \right\|_{H^2(M)}\underset{t\rightarrow +\infty}{\longrightarrow} 0$$ and $v_\infty \in C^{\infty}(M)$ is a solution of 
\begin{equation}
\label{nonex1}
-\Delta v_\infty +Q=\frac{\rho e^{v_\infty}}{\int_M e^{v_\infty}dV}.
\end{equation}
It follows that 
\begin{equation*}
\left\|v_\infty \right\|_{C^{2+\alpha}(M)}\leq C_0,
\end{equation*} 
where $C_0$ is the constant defined in \eqref{12nov2012e1}.
This implies that there exists a constant $\bar{C}$ depending on $M,Q,\rho$ and $C_0$ such that
$$J(w_\infty)\geq -\bar{C}.$$
Since $J_\rho (v(t_1))\leq J_\rho (v(t_2))$, for all $t_1 \geq t_2$, we have
$$J_\rho(v_0)\geq J_\rho (v_\infty )\geq -\bar{C}.$$
But, by hypothesis, $J_\rho (v_0)\leq -C$, choosing $C >\bar{C}$, we get a contradiction.  

\bibliographystyle{plain}
\bibliography{biblio}
\end{document}